\documentclass[11pt]{article}
\usepackage{fullpage}

\usepackage{amsmath,amsthm,amssymb}
\usepackage[algo2e,ruled]{algorithm2e}

\newtheorem{theorem}{Theorem}
\newtheorem{corollary}{Corollary}
\newtheorem{lemma}{Lemma}
\newtheorem{proposition}{Proposition}

\newtheorem{assumption}{Assumption}

\newcommand{\ip}[1] {\langle #1 \rangle }

\newcommand{\E}{\mathbf{E}}
\newcommand{\R}{\mathbf{R}}

\newcommand{\cA}{\mathcal{A}}

\newcommand{\cX}{\mathcal{X}}

\newcommand{\cS}{\mathcal{S}} 
    
\newcommand{\cO}{{\mathcal{O}}}

\newcommand{\cY}{\mathcal{Y}}

\newcommand{\argmin}{\operatornamewithlimits{argmin}}
\newcommand{\argmax}{\operatornamewithlimits{argmax}}

\newcommand{\normm}[1]{{\vert\kern-0.25ex\vert\kern-0.25ex\vert #1 \vert\kern-0.25ex\vert\kern-0.25ex\vert}}
\newcommand{\normb}[1]{{\big\vert\kern-0.25ex\big\vert\kern-0.25ex\big\vert #1 \big\vert\kern-0.25ex\big\vert\kern-0.25ex\big\vert}}
\newcommand{\normbb}[1]{{\bigg\vert\kern-0.25ex\bigg\vert\kern-0.25ex\bigg\vert #1 \bigg\vert\kern-0.25ex\bigg\vert\kern-0.25ex\bigg\vert}}

\newcommand{\xs}{ \hat{x} }
\newcommand{\ys}{ \hat{y} }
\newcommand{\xsi}{ \hat{x}_{(i)} }
\newcommand{\ysi}{ \hat{y}_{(i)} }

\title{Generalization Bounds for Stochastic Saddle Point Problems}

\author{%
	Junyu Zhang\thanks{University of Minnesota, Twin Cities, Minneapolis, Minnesota 55455(\texttt{zhan4393@umn.edu})}
	\and
	Mingyi Hong \thanks{University of Minnesota, Twin Cities, Minneapolis, Minnesota 55455(\texttt{mhong@umn.edu})}
	\and
	Mengdi Wang\thanks{Princeton University,  Princetone, NJ 08544 (\texttt{mengdw@princeton.edu})}
	\and
	Shuzhong Zhang \thanks{University of Minnesota, Twin Cities,
		Minneapolis, Minnesota 55455(\texttt{zhangs@umn.edu})}
}

\date{\today} 

\begin{document}
\maketitle 

\begin{abstract}
	This paper studies the generalization bounds for the empirical saddle point (ESP) solution to stochastic saddle point (SSP) problems. For SSP with Lipschitz continuous and strongly convex-strongly concave objective functions, we establish an {\footnotesize$\cO\left(1/n\right)$} generalization bound by using a uniform stability argument. We also provide generalization bounds under a variety of assumptions, including the cases without strong convexity and without bounded domains. We illustrate our results in two examples: batch policy learning in Markov decision process, and mixed strategy Nash equilibrium estimation for stochastic games. In each of these examples, we show that a regularized ESP solution enjoys a near-optimal sample complexity. 
	To the best of our knowledge, this is the first set of results on the generalization theory of ESP. 
\end{abstract}
 
\section{Introduction}
\label{sec:intro}
Consider the stochastic saddle point (SSP) problem
\begin{equation}
\label{prob:SSP}
(\text{SSP})\qquad\qquad\qquad\qquad\qquad\min_{x\in\cX}\max_{y\in\cY} \Phi(x,y):=\E\big[\Phi_\xi(x,y)\big] ,\qquad\qquad\qquad\qquad\qquad
\end{equation}
where $\cX$ and $\cY$ are compact and convex sets, and $\xi$ is a random variable. We denote the optimal solution to \eqref{prob:SSP} as $(x^*,y^*)$.
SSP problem finds a wide range of applications in machine learning, reinforcement learning, operations research and game theory. Many stochastic approximation (SA) algorithms have been proposed for approximating the SSP solution based on samples of $\xi$, see e.g. \cite{natole2018stochastic,nemirovski2009robust,shalev2013stochastic,xiao2019dscovr,yan2019stochastic,zhang2017stochastic,zhao2019optimal}. 
Most of the algorithms make primal-dual updates and enjoy convergence guarantees with appropriately chosen stepsizes.

In this paper, we provide understanding of \eqref{prob:SSP} from a finite sample perspective. Consider the empirical counterpart of the SSP problem, which we refer to as the empirical saddle point (ESP) problem
\begin{equation}
\label{prob:ESP}
(\text{ESP})\qquad\,\qquad\quad\qquad\qquad\min_{x\in\cX}\max_{y\in\cY} \hat\Phi_n(x,y):=\frac{1}{n}\sum_{i=1}^n\Phi_{\xi_i}(x,y) ~,\qquad\qquad\qquad\qquad\quad
\end{equation}
where $\Gamma:= \{\xi_1,...,\xi_n\}$ is a collection of $n$ i.i.d. samples of $\xi$'s. A natural way of estimating the optimal solution $(x^*,y^*)$ to \eqref{prob:SSP} is to solve instead its empirical approximation given by \eqref{prob:ESP}. This approach is also known as sample average approximation (SAA), see \cite{shapiro2014lectures}. We denote by $(\hat x,\hat y)$ the empirical saddle point (ESP) solution to problem \eqref{prob:ESP}. Based on a given set of samples, one can compute the ESP solution using any convergent algorithm for minimax optimization.  In parallel to the generalization theory for empirical risk minimization \cite{vapnik1992principles,vapnik2006estimation,vapnik2013nature}, we aim to analyze the empirical saddle point and establish finite-sample generalization bounds.

\subsection{Motivating Examples}

Stochastic saddle points \eqref{prob:SSP} are very common in machine learning, game theory, and operations research. A generalization theory for ESP would be useful for establishing generalization bounds for a number of machine learning tasks that are beyond empirical risk minimization. We will study two examples in this paper.

One example is batch policy learning for Markov Decision Process (MDP). For the infinite-horizon average-reward MDP, the policy optimization problem is equivalent to an SSP problem, known as the Bellman saddle point \cite{puterman2014markov,wang2017primal}, given by 
\begin{eqnarray}
\label{prob:aMDP}
\min_{x\in\cX}\max_{y\in\cY} \Phi(x,y): =\ip{y,r}+ \sum_{a\in\cA}y_a^\top(P_a-I)x,
\end{eqnarray} 
where $x$ is the value function, $y=\{y_a\}_{a\in\cA}$ is the state-action occupancy measure, $a$ and $\cA$ denote the action and action space respectively, $P_a$ denotes the transition probability matrix under $a$, $r$ is the reward function (see Section \ref{sec:aMDP} for details). In the batch policy learning problem, we want to solve the MDP without knowledge of $r,P_a$, instead we only have sample state transitions. This motivates us to study the empirical optimal policy that is equivalent to the ESP solution of the Bellman saddle point \eqref{prob:aMDP}. 

Another example is   the two-person stochastic matrix game
\begin{equation}
\label{prob:game}
\min_{x\in\Delta_{N_1}}\max_{y\in\Delta_{N_2}} x^\top\E[A_\xi]y\qquad\mbox{with}\qquad\Delta_{N_i}:=\{z\in\R^{N_i}: z\geq0, \mathbf{1}^\top z=1\},
\end{equation}
where $x,y$ are the mixed strategies of the two players, and $A_\xi\in\R^{N_1\times N_2}$ is a stochastic payoff matrix. Based on sample payoffs from past plays, we can estimate the mixed strategy Nash equilibrium by solving a regularized version of the empirical matrix game. Our generalization bound will be used to evaluate the quality of the empirical Nash equilibrium learned from data. 

Additional examples include the stochastic composite optimization, which finds applications in off-policy policy evaluation and risk-averse optimization, of the form
\begin{equation*}
\label{prob:sto-compOpt}
\min_{x\in\cX} \,\, r(x) + f\left(\E\big[A_\xi x - b_\xi\big]\right),
\end{equation*}
where $A_\xi$ is a random matrix, $b_\xi$ is a random vector, $f$ is a convex loss; see \cite{Composite-Nesterov,Composite-Wang,Composite-zhang} and references therein. By using the convex conjugate $f^*$ of $f$, we can reformulate the problem as an SSP problem:
\begin{equation*}
\min_{x\in\cX}\max_{y\in\cY} r(x) + \E\left[\big(A_\xi x-b_\xi\big)^\top y\right] - f^*(y).
\end{equation*}
Thus a generalization theory for this SSP problem would lead to generalization bounds for the original composite optimization problem. Similarly, for the robust optimization problem, where the aim is to minimize the worst loss among multiple scenarios,
\begin{equation*}
\min _{x\in\cX} F(x):= \max\left\{\E_\xi\big[f_{1,\xi}(x)\big],\E_\xi\big[f_{2,\xi}(x)\big],\cdots,\E_\xi\big[f_{m,\xi}(x)\big]\right\},
\end{equation*}
which can be reformulated as the following SSP problem:
\begin{equation*}
\min_{x\in\cX}\max_{y\geq0, \mathbf{1}^\top y = 1} \sum_{i=1}^m \E_\xi\big[y_i\cdot f_{i,\xi}(x)\big],
\end{equation*}
where we have applied the fact that  $\max_i \{z_i\} = \max\left\{y^\top z: y\geq0, \mathbf{1}^\top y = 1\right\}$.

\subsection{Weak and Strong Generalization Measures}

We will study the generalization properties of the empirical saddle point (ESP) solution $(\xs,\ys)$  via two metrics.  The first metric, referred to as the \emph{weak generalization measure \footnote{We decide not to use ``weak (strong) duality gap measure'' to avoid confusion with the well known terminologies of weak (strong) duality.} (WGM)}, is defined as 
\begin{equation}
\label{defn:measure-weak}
\Delta^{w}(\xs,\ys): = \max_{y\in\cY}\E\big[\Phi(\xs,y)\big] - \min_{x\in\cX} \E\big[\Phi(x,\ys)\big],
\end{equation}
where the expectations are taken over the sample set $\{\xi_1,\ldots,\xi_n\}$.
In some applications, one desires a stronger metric of optimality, which we refer to as the \emph{strong generalization measure (SGM)}, given by
\begin{equation}
\label{defn:measure-strong}
\Delta^s(\xs,\ys):=\E\Big[\max_{y\in\cY}\Phi(\xs,y) - \min_{x\in\cX} \Phi(x,\ys)\Big].
\end{equation}
The SGM is usually referred to as the expected \emph{duality gap} in the optimization community. A third commonly used metric is $d^2(\xs,\ys):=\| \xs -x^*\|^2+\normm{\ys-y^*}^2$, which is the squared distance between the ESP solution to the true saddle point solution. Note that here we allow the use of two different norms $\|\cdot\|$ and $\normm{\cdot}$ to measure the distances in $\cX$ and $\cY$, respectively. Suppose that $\Phi(\cdot,\cdot)$ is strongly convex-strongly concave (SC-SC). Then we have $\Delta^w(\xs,\ys)\geq \Omega\big(\E[d^2(\xs,\ys)]\big)$. Due to the Jensen's inequality, we also have $\Delta^w(\xs,\ys)\leq\Delta^s(\xs,\ys)$. Therefore the SGM is strongest among the three, WGM is the second and $d^2(\xs,\ys)$ is the weakest, i.e.,
$$\Omega\left(\E[d^2(\xs,\ys)]\right)\leq \Delta^w(\xs,\ys) \leq  \Delta^s(\xs,\ys).$$
Based on the goal of the specific problem studied, different generalization measure may be desired. When it suffices to bound the WGM, e.g. the MDP \eqref{prob:aMDP}, analyzing the problem by bounding SGM will not result in a tight sample complexity.  

\subsection{Main Results}

\begin{table}	
	\centering
	\caption{Generalization bounds for stochastic saddle points (results from this paper) \label{table:summary}}
	\begin{tabular}{c | c | c | c}
		\hline
		& Assumption \ref{assumption:convex-concave} \& \ref{assumption:Lipschitz}  & Assumption \ref{assumption:convex-concave}, \ref{assumption:Lipschitz} \& \ref{assumption:grad-Lipschitz}  &Assumption \ref{assumption:convex-concave}, \ref{assumption:grad-Lipschitz} \& \ref{assumption:grad-opt}  \\ \hline\hline
		$\E\left[\!d^2(\xs,\ys)\!\right]$  & $\cO\left(\frac{1}{n\mu_{\min}}\bigg(\frac{(\ell_x^w)^2}{\mu_x}+\frac{(\ell_y^w)^2}{\mu_y}\bigg)\right)$  & Same as left   & $\cO\left(\frac{\kappa^2\cdot\E\left[\!\!\|\!\nabla\!\Phi_{\!\xi}\!(x^*\!\!,y^*\!)\!\|^2\!\right]}{n\cdot\min\{\mu_x^2,\mu_y^2\}}\right)$   \\ \hline
		$\Delta^w\!\left(\xs,\ys\right)$ \!\!\!\!  & $\cO\left(\frac{(\ell_x^w)^2}{n\mu_x}+\frac{(\ell_y^w)^2}{n\mu_y}\right)$   & Same as left &  $\cO\left(\frac{\kappa^4\cdot\E\left[\!\!\|\!\nabla\!\Phi_{\!\xi}\!(x^*\!\!,y^*\!)\!\|^2\!\right]}{n\cdot\min\{\mu_x,\mu_y\}}\right)$   \\ \hline
		$\Delta^s\!\left(\xs,\ys\right)$ \!\!\! \! & NA & $\cO\left(\!\sqrt{\frac{L^2_{xy}}{\mu_x\mu_y}\!+\!1}\bigg(\frac{(\ell_x^s)^2}{n\mu_x}\!+\!\frac{(\ell_y^s)^2}{n\mu_y}\bigg)\!\right)$ &  $\cO\left(\frac{\kappa^4\cdot\E\left[\!\!\|\!\nabla\!\Phi_{\!\xi}\!(x^*\!\!,y^*\!)\!\|^2\!\right]}{n\cdot\min\{\mu_x,\mu_y\}}\right)$  \\ \hline
	\end{tabular}
\end{table}

In this paper, we establish the generalization bounds for solving SSP problem using the empirical saddle point solution under various assumptions. See Table \ref{table:summary} for an overview of our technical results. Contributions of the paper are three-folded:
\begin{itemize}
	\item We establish a uniform stability argument for the ESP solution by extending the technique of \cite{shalev2009stochastic}. For SSP problem over compact domains that are Lipschitz continuous and SC-SC, we provide an $\cO(1/n)$ bound for the WGM metric.  With an additional assumption on gradient Lipschitz continuity, we also establish an $\cO(1/n)$ bound for the SGM metric. 
	\item Further, we extend the generalization theory to SSP problems with unbounded domains or without the SC-SC property. By using a different analysis, we establish an $\cO(1/n)$ generalization bound for the ESP solution even if the feasible regions are unbounded. We also provide a generalization bound for the regularized ESP problem. 
	\item We illustrate the applications of our theory in two important learning tasks: the batch policy learning in MDP, and the Nash equilibrium estimation. In each of these tasks, the empirical saddle point provides a conceptually simple estimator. Further, we use the generalization theory to show these estimators have provable sample complexities of
	$\tilde{\cO}\left(\frac{|\cS||\cA|}{\epsilon^2}\right)$ and  $\tilde{\cO}\left(\frac{N_1N_2}{\epsilon^{2}}\right)$ for the two learning tasks, respectively\footnote{The definition of these parameters can be found in the corresponding sections.}. These sample complexities have tight dependences on the problems' parameters, and they have not been studied before. 
\end{itemize}

\subsection{Related works}
Let us review the stochastic approximation (SA) approaches for the SSP problem \eqref{prob:SSP}. When $\Phi(\cdot,\cdot)$ is only convex and concave, Nemirovski et al. \cite{nemirovski2009robust} established an $\cO(1/\sqrt{n})$ convergence in SGM for a  stochastic mirror descent ascent algorithm. Similar $\cO(1/\sqrt{n})$ convergence in SGM are also obtained by Bach \& Levy \cite{bach2019universal} and Zhao et al. \cite{zhao2019optimal}, under various assumptions. When the SC-SC property is further assumed, a faster $\cO(1/n)$ convergence can be derived. For example,  Natole et al. \cite{natole2018stochastic} obtained an $\cO(1/n)$ convergence in terms of the squared distance metric. Yan et al. \cite{yan2019stochastic} designed a stochastic gradient method with $\cO(1/n)$ convergence in SGM when the coupling between the primal and dual variable is linear. Yan et al. \cite{yan2020sharp} derived an epoch-wise stochastic gradient method that guarantees $\cO(1/n)$ in SGM. This result of \cite{yan2020sharp} does not rely on the linear coupling structure, but instead requires additional conditions on the tail distribution of sampling noise. There also exist research results for SSP problem with special structures such as finite-sum and  bilinear coupling, etc, see e.g. \cite{du2018linear,shalev2013stochastic,xiao2019dscovr,zhang2017stochastic}.

There exits a rich body of literatures on the generalization theory for solving stochastic convex optimization (SCO) by empirical risk minimization (ERM), namely, 
$$\min_{x\in\cX} \Phi(x):=\E_\xi\big[\Phi_{\xi}(x)\big]\quad \mbox{(SCO)}\qquad\quad\mbox{ and } \quad\qquad\min_{x\in\cX} \hat\Phi_n(x):=\frac{1}{|\Gamma|}\sum_{\xi\in\Gamma}\Phi_{\xi}(x)\quad\mbox{(ERM)}.$$
In the seminal paper \cite{shalev2009stochastic}, Shalev et al. established an $\cO(1/n)$ ERM generalization bound for strongly convex problems and an $\cO(1/\sqrt{n})$ risk bound for general convex problems. Similar rates are also obtained in related works \cite{sridharan2009fast,gonen2017average}. The main technique used by \cite{shalev2009stochastic} and our paper is the \emph{uniform stability} argument, which was originally introduced by Bousquet and Elisseef \cite{bousquet2002stability}, and later on studied in many papers, see e.g. \cite{kearns1999algorithmic,mukherjee2006learning,shalev2009stochastic,shalev2013stochastic,hardt2015train,chen2018stability}, etc.  With the tool of the Rademacher complexity $R_n$, Srebro et al. \cite{srebro2010smoothness} demonstrated an $\cO(R_n/\sqrt{n})$ risk bound for ERM, and many papers strengthened the theory further \cite{bartlett2002rademacher,bartlett2005local}. For nonconvex but exp-concave objectives, Koren \& Levy \cite{koren2015fast} and Mehta \cite{mehta2016fast} derived a risk bound of $\cO(1/n)$.  Under certain stronger conditions, a tighter $\cO(1/n^2)$ risk bound has been shown \cite{zhang2017empirical}. To the authors' best knowledge, there is no existing generalization bound for stochastic saddle point problems.


\section{Generalization Bounds for Empirical Saddle Points}
\label{sec:main}
\subsection{Assumptions}
In most of our analysis, we require that the objective function $\Phi$ is strongly convex and strongly concave (SC-SC), as stated in the following assumption.
\begin{assumption}[SC-SC objective function]
	\label{assumption:convex-concave}
	$\exists\mu_x,\mu_y>0$, s.t. for almost every $\xi$, $\Phi_{\!\xi}(\cdot,y)$ is $\mu_x$-strongly convex under norm $\|\cdot\|$ and $\Phi_\xi(x,\cdot)$ is $\mu_y$ strongly concave under the norm $\normm{\cdot}$. Namely,
	\begin{equation}
	\begin{cases}
	\Phi_\xi(x',y)\geq\Phi_\xi(x,y)  + \langle u,x'-x\rangle + \frac{\mu_x}{2}\|x'-x\|^2, & \forall x,x'\in\cX ,y\in\cY\,\,\mbox{and}\,\,u\in\partial_x\Phi_\xi(x,y),\\
	\Phi_\xi(x,y')\leq\Phi_\xi(x,y) + \langle v,y'-y\rangle - \frac{\mu_y}{2}\normm{y'-y}^2,& \forall y,y'\in\cY,x\in\cX\,\,\mbox{and}\,\,v\in\partial_y \Phi_\xi(x,y),
	\end{cases}
	\end{equation} 
	where $\partial_x\Phi_\xi(\cdot,y)$ and $\partial_y\Phi_{\!\xi}(x,\cdot)$ denote the subgradients and supergradients, respecitively. 
\end{assumption}
For convex analysis of strongly convex function under non-Euclidean norms, see \cite{shalev2007online,kakade2012regularization} and references therein. We further assume that the feasible regions $\cX, \cY$ are bounded and the objective function is  Lipschitz continuous. 
\begin{assumption}[Function Lipschitz continuity]
	\label{assumption:Lipschitz}
	The feasible regions $\cX$ and $\cY$ are compact convex sets. For almost every $\xi$,  there exist constants $\ell_x(\xi,y)$ and $\ell_y(\xi,x)$ s.t. 
	\begin{equation}
	\begin{cases}
	|\Phi_\xi(x',y)-\Phi_\xi(x,y)| \leq \ell_x(\xi,y)\|x'-x\|,& \forall x,x'\in\cX \,\,\mbox{and}\,\,y\in\cY,\\
	|\Phi_\xi(x,y')-\Phi_\xi(x,y)|\leq \ell_y(\xi,x)\normm{y'-y},& \forall y,y'\in\cY\,\,\mbox{and}\,\,x\in\cX.
	\end{cases}
	\end{equation} 
	To bound the WGM, we assume
	\begin{equation}
	\label{assumption-weak}
	(\ell_x^w)^2:=\sup_{y\in\cY}\E_\xi\big[\ell^2_x(\xi,y)\big]<+\infty\quad\mbox{and}\quad(\ell_y^w)^2:=\sup_{x\in\cX}\E_\xi\big[\ell^2_y(\xi,x)\big]<+\infty.
	\end{equation}
	To bound the SGM, we assume
	\begin{equation}
	\label{assumption-strong}
	(\ell_x^s)^2:=\E_\xi\big[\sup_{y\in\cY}\ell^2_x(\xi,y)\big]<+\infty\quad\mbox{and}\quad(\ell_y^s)^2:=\E_\xi\big[\sup_{x\in\cX}\ell^2_y(\xi,x)\big]<+\infty.
	\end{equation}
	Due to Jensen's inequality, $\ell_x^w\leq \ell_x^s$ and $\ell_y^w\leq\ell_y^s$ always hold.
\end{assumption}

In our analysis, Assumptions \ref{assumption:convex-concave} and \ref{assumption:Lipschitz} only guarantee the $\cO(1/n)$ bound for the WGM metric. In order to prove an $\cO(1/n)$ bound for the stronger metric SGM, we will require additional smoothness of $\Phi$.
\begin{assumption}[Gradient Lipschitz continuity]
	\label{assumption:grad-Lipschitz}
	There exist constants $L_x$, $L_y$ and $L_{xy}$ s.t. for $\forall x,x',y,y'$, it holds that 
	\begin{eqnarray*}
		&&\!\!\!\!\!\!\!\!\|\nabla_x \Phi(x,y) - \nabla_x \Phi(x',y) \|_*\leq L_x\|x-x'\|,\qquad\,
		\normm{\nabla_y \Phi(x,y) - \nabla_y \Phi(x,y')}_*\leq L_y\normm{y-y'},\qquad\\
		&&\!\!\!\!\!\!\!\!\|\nabla_x \Phi(x,y) - \nabla_x \Phi(x,y') \|_*\leq L_{xy}\normm{y-y'},\qquad\!\!\!
		\normm{\nabla_y \Phi(x,y) - \nabla_y \Phi(x',y) }_*\leq L_{xy}\|x-x'\|,
	\end{eqnarray*}
	where $\|\cdot\|_*$ and $\normm{\cdot}_*$ stand for the dual norms of $\|\cdot\|$ and $\normm{\cdot}$ respectively.
\end{assumption}
Finally, we also study the case where $\cX$ and $\cY$ are unbounded. Such unboundedness would invalidate Assumption \ref{assumption:Lipschitz} as well as the stability argument. To remedy this issue, we would replace Assumption \ref{assumption:Lipschitz} with the following assumption about the true optimal solution $(x^*,y^*)$.

\begin{assumption}
	\label{assumption:grad-opt}
	There exists a constant $C$ s.t. $\E_\xi\left[\|\nabla \Phi_\xi(x^*,y^*)\|_2^2\right]\leq C<+\infty.$
\end{assumption}

\subsection{Main Results}
We use the leave-one-out technique in \cite{shalev2009stochastic} to analyze the stability of the ESP solutions. Let $\Gamma:=\{\xi_1,...,\xi_n\}$ be a set of $n$ i.i.d. samples, and let $\xi_i'$ be another independent sample. We then define the perturbed sample set $\Gamma(i) = \Gamma\cup\{\xi_i'\}\backslash\{\xi_i\}.$ That is, 
$\Gamma(i)$ is constructed by replacing just the $i$-th sample $\Gamma$ with another i.i.d. sample $\xi_i'$.  In the next lemma, we establish the stability property of the ESP solutions.
\begin{lemma}[Stability property]
	\label{lemma:stability-SC-SC}
	Let the Assumptions \ref{assumption:convex-concave} and  \ref{assumption:Lipschitz} hold. Denote $(\xs,\ys)$ and $(\xsi,\ysi)$ as the solutions to the ESP problems with sample sets $\Gamma$ and $\Gamma(i)$ respectively. Then
	$$
	\sqrt{\mu_x\|\xs-\xsi\|^2 + \mu_y\normb{\ys-\ysi}^2} 
	\leq  \frac{1}{n} \sqrt{\frac{(\ell_x(\xi_i,\ysi) + \ell_x(\xi_i',\ys))^2}{\mu_x} + \frac{(\ell_y(\xi_i,\xsi) + \ell_y(\xi_i',\xs))^2 }{\mu_y} }.
	$$
\end{lemma}
Note that Lemma \ref{lemma:stability-SC-SC} can be reduced from Lemma \ref{lemma:stability-reg} by setting the regularizer to be 0, hence we only provide the proof of the latter result, see Appendix \ref{appdx:stability-reg}. Note that by the Mcdiarmid's inequality \cite{McD1,McD2}, the stability argument of Lemma \ref{lemma:stability-SC-SC} immediately results in an $\tilde{\cO}(1/\sqrt{n})$ generalization bound, which, however, is not tight. In Theorem \ref{theorem:gen-SC-SC-weak}, we establish a tighter $\cO(1/n)$ bound by using a more careful analysis. See Appdendix \ref{appdx:thm-gen-sc-sc-weak} for the proof. 
\begin{theorem}(Upper bound on WGM)
	\label{theorem:gen-SC-SC-weak}
	For the SSP problem \eqref{prob:SSP}, let $(\xs,\ys)$ be the solution to the ESP problem \eqref{prob:ESP}. And let $n = |\Gamma|$ be sample size. Under Assumptions \ref{assumption:convex-concave} and \ref{assumption:Lipschitz}, we have 
	\begin{equation}
	\label{thm:gen-weak}
	\E \left[d^2(\xs,\ys)\right]\leq \frac{\Delta^w(\xs,\ys)}{\min\{\mu_x,\mu_y\}}\qquad\mbox{ and }\qquad \Delta^w(\xs,\ys) \leq \frac{2\sqrt{2}}{n}\cdot\left( \frac{(\ell^w_x)^2}{\mu_x} + \frac{(\ell_y^w)^2}{\mu_y}\right).
	\end{equation}
\end{theorem}

The generalization bounds given in Theorem \ref{theorem:gen-SC-SC-weak} have tight dependence on the sample size $n$, as well as the problem's parameters $\ell_x^w, \ell_y^w$ and $\mu_x, \mu_y$. To see this, we can simply consider the special case of SCO.
When $\Phi_{\xi}(x,\cdot)\equiv f_\xi(x)$, i.e., the objective function is constant in $y$, the SSP and ESP reduce to the classical SCO and ERM respectively. In this case, the difference between $\ell_x^w$ and $\ell_x^s$ vanishes, and we denote them as $\ell_x: = \ell_x^w = \ell_x^s$. The WGM also reduces to
$$\Delta^w(\xs,\ys) = \max_{y\in\cY}\E[f(\xs)] - \min_{x\in\cX}\E[f(x)] = \E\big[f(\xs) - \min_{x\in\cX}f(x)\big].$$
It follows that the generalization bound \eqref{thm:gen-weak} becomes $\cO\left(\frac{\ell_x^2}{n\mu_x}\right)$ and matches the known generalization lower bound for ERM \cite{shalev2009stochastic}.

By utilizing additional smoothness of the objective function, we provide an $\cO(1/n)$ bound on SGM in the following theorem, whose proof is given in Appendix \ref{appdx:thm-gen-sc-sc-strong}.
\begin{theorem}(Upper bound on SGM)
	\label{theorem:gen-SC-SC-strong}
	Under the settings of Theorem \ref{theorem:gen-SC-SC-weak}, if Assumption \ref{assumption:grad-Lipschitz} holds in addition, we have 
	\begin{equation}
	\label{thm:gen-strong}
	\Delta^s(\xs,\ys) \leq \frac{2\sqrt{2}}{n}\cdot\sqrt{\frac{L^2_{xy}}{\mu_x\mu_y}+1}\cdot\left(\frac{(\ell_x^s)^2}{\mu_x}+\frac{(\ell_y^s)^2}{\mu_y}\right).
	\end{equation}
\end{theorem}
We remark that the bound \eqref{thm:gen-strong} has an additional multiplicative $\cO\big(L_{xy}/\sqrt{\mu_x\mu_y}\big)$ factor compared to bound \eqref{thm:gen-weak}. It remains open whether this dependence can be improved, as a question for future work. 

It is worth mentioning that both Theorems \ref{theorem:gen-SC-SC-weak} and \ref{theorem:gen-SC-SC-strong} are based on the stability argument in Lemma \ref{lemma:stability-SC-SC}, which relies heavily on the  Lipschitz continuity of the objective function. Next we study SSP problems over unbounded domains. In this case, the SC-SC property and function Lipschitz continuity are mutually exclusive. In the next theorem, we study the SSP problems with unbounded domains and provide a generalization bound without assuming Lipschitz continuity of the objective.

\begin{theorem}[Generalization error for unbounded problems]
	\label{theorem:gen-unbounded}
	Let Assumptions \ref{assumption:convex-concave}, \ref{assumption:grad-Lipschitz} and \ref{assumption:grad-opt} hold, and let the SSP be unconstrained. Let $\|\cdot\| = \normm{\cdot} = \|\cdot\|_2$. Then 
	\begin{eqnarray}
	\E\left[d^2\left(\xs,\ys\right)\right]
	\leq \cO\left(\frac{C\kappa^2}{n\mu^2}\right)\qquad\mbox{and}\qquad\Delta^s(\xs,\ys)
	\leq  \cO\left(\frac{C\kappa^4}{n\mu}\right)
	\end{eqnarray}	
	where $\mu = \min\{\mu_x, \mu_y\}$ and $\kappa = \frac{\max\{L_x, L_y, L_{xy}\}}{\min\{\mu_x,\mu_y\}}$  is the condition number. 
\end{theorem} 

\subsection{Generalization Bounds for Regularized ESP}

Next we study SSP problems that are not necessarily strongly convex or strongly concave (SC-SC). 
In order to get a stable solution, we consider a regularized version of the empirical saddle point problem
\begin{equation}
\label{prob:ESP-reg}
(\text{R-ESP})\qquad\quad\qquad\quad\,\qquad\qquad\min_{x\in\cX}\max_{y\in\cY} \hat\Phi_n(x,y) + \Psi(x,y),\qquad\qquad\qquad\quad\qquad\quad\quad
\end{equation}
where $\Psi$ is a regularization function that is SC-SC and can be specified by the user.  
Let us generalize the results of Lemma \ref{lemma:stability-SC-SC} and Theorem \ref{theorem:gen-SC-SC-weak} to the case with regularization.

\begin{lemma}[Stability property for regularized ESP solution]
	\label{lemma:stability-reg}
	Suppose the regularization function $\Psi(\cdot,y)$ is $\nu_x$-strongly convex under norm $\|\cdot\|$ and $\Psi(x,\cdot)$ is $\nu_y$-strongly concave under norm $\normm{\cdot}$, and there exists $R>0$ s.t. $|\Psi(x,y)|\leq R,$  $\forall (x,y)\in\cX\times\cY.$
	Under the Assumptions \ref{assumption:convex-concave},  \ref{assumption:Lipschitz}, let $(\xs,\ys)$ and $(\xsi,\ysi)$ be the solution to the R-ESP problem \eqref{prob:ESP-reg} with sample set $\Gamma$ and $\Gamma(i)$ respectively, then
	\begin{eqnarray*}
		\sqrt{(\mu_x\!+\!\nu_x)\|\xs-\xsi\|^2 + (\mu_y\!+\!\nu_y)\normm{\ys-\ysi}^2}
		\leq \frac{1}{n}\sqrt{\frac{(\ell_x(\xi_i,\!\ysi) \!+\! \ell_x(\xi_i',\!\ys))^2\!}{\mu_x\!+\!\nu_x} \!+\! \frac{(\ell_y(\xi_i,\!\xsi) \!+\! \ell_y(\xi_i',\!\xs))^2 }{\mu_y+\nu_y}}.
	\end{eqnarray*}
\end{lemma} 
Note that Lemma \ref{lemma:stability-reg} is not a trivial extension of Lemma \ref{lemma:stability-SC-SC}: The Lipschitz constant of $\Psi$ does not contribute to the stability bound in Lemma \ref{lemma:stability-reg}. Its proof is given in \ref{appdx:stability-reg}. Then we obtain the following generalization bound for R-ESP, whose proof is given in Appendix \ref{appdx:weak-gen-reg}.
\begin{lemma}[Generalization bound for regularized ESP]
	\label{theorem:weak-gen-reg}
	Under the settings of Lemma \ref{lemma:stability-reg}, the R-ESP solution $(\hat x,\hat y)$ satisfies 	
	$$\Delta^w(\xs,\ys) \leq \frac{2\big(\ell_x^w\big)^2}{n(\mu_x+\nu_x)} + \frac{2\big(\ell_y^w\big)^2}{n(\mu_y+\nu_y)} + 2R,$$
	where the WGM $\Delta^w(\cdot,\cdot)$ is defined for the original unregularized SSP problem.
\end{lemma}

Finally we show how to choose the regularizer optimally and establish a generalization bound that depends only on the diameters of $\cX,\cY$ and constants of Lipschitz continuity of $\Phi$.
\begin{corollary} Suppose the $\Phi_{\xi}$'s are convex concave but not SC-SC. Namely, $\mu_x = \mu_y = 0$. Suppose the $\Phi_\xi$'s satisfy Assumption \ref{assumption:Lipschitz} under the $L_2$ norm. Then we can set the regularizer to be 
	$\Psi(x,y) = \frac{\alpha_x}{2}\|x\|_2^2 - \frac{\alpha_y}{2}\|y\|_2^2.$
	Consequently, $R = \frac{\alpha_x}{2}D_x^2 + \frac{\alpha_y}{2}D_y^2$, where $D_x$ and $D_y$ denote the diameters of $\cX$ and $\cY$ under $L_2$-norm respectively. If we set $\alpha_x = \frac{\ell_x^w}{\sqrt{n}D_x}$ and $\alpha_y = \frac{\ell_y^w}{\sqrt{n}D_y}$, Lemma \ref{theorem:weak-gen-reg}  implies
	$$\Delta^w(\xs,\ys) \leq \cO\left(\frac{\ell_x^wD_x + \ell_y^wD_y}{\sqrt{n}}\right).$$
\end{corollary}


\section{Application to Batch Policy Learning for MDP}

\label{sec:aMDP}
\subsection{Saddle Point Formulation of MDP}
Consider the policy learning problem for an infinite-horizon average-reward Markov Decision Process (MDP). The MDP instance is specified by $M = (\cS,\cA,P, r)$ where $\cS$ is a finite state space. $\cA$ is a finite action space, $P=\{P_a\}$ are state transition matrices with $P_a(s,s') = \mathbf{Prob}(s_{t+1} = s'\,\,|\,\,s_t = s, a_t = a)$, for $\forall s,s'\in\cS$ and $a\in\cA$. $r$ is the reward function with $r_{sa}\in[0,1]$ being the reward received after taking action $a$ at state $s$. A policy $\pi:\cS\mapsto \Delta_{\cA}$ maps a state $s$ to a distribution over the action space $\cA$, where we denote the probability of taking action $a$ at state $s$ as $\pi(a|s)$. The objective is to maximize the long-term average reward, defined as
\begin{equation}
\label{defn:value-func-opt-1}
\hat v^*:=\max_{\pi} \lim_{T\rightarrow\infty}\E\left[\frac{1}{T}\sum_{t=0}^{T-1} r_{s_t a_t}~\bigg|~\pi, s_0 = s\right].
\end{equation}
The optimal Bellman equation has an equivalent saddle point formulation \eqref{prob:aMDP} \cite{puterman2014markov,chen2016stochastic}
\begin{equation*} 
\min_{x\in\cX}\max_{y\in\cY} \Phi(x,y): =\ip{y,r}+\sum_{a\in\cA} y_a^\top(P_a-I)x,
\end{equation*}
where $x\in\R^{|\cS|}$ is the {\it difference-of-value} vector, $y\in\R^{|\cS|\times|\cA|}$ stands for the stationary state-action distribution under certain policy $\pi$. $y_a = [y_{1,a},...,y_{|\cS|,a}]'$ is the $a$-th column of $y$. Under the assumption of fast mixing time and uniform ergodicity (Assumptions 2,3 of \cite{chen2018scalable}), there exists constant $t_{mix}$ and $\tau$ such that one can set the feasible regions $\cX$ and $\cY$ as
\begin{equation}
\label{defn:sets-aMDP} 
\begin{cases}
\cX:=\{x\in\R^{|\cS|}: \|x\|_\infty\leq 2t_{mix}\},\\
\cY := \big\{y\in\R^{|\cS|\times|\cA|}\,:\,y\geq 0,\, \|y\|_1=1,\, \frac{1}{\sqrt{\tau}|\cS|}\cdot\mathbf{1}\leq \sum_{a\in\cA}y_a\leq \frac{\sqrt{\tau}}{|\cS|}\cdot\mathbf{1}\big\},
\end{cases}
\end{equation}
(see Appendix \ref{appdx:aMDP-assumptions} for details). In the policy learning setting, we do not know either $P$ or $r$. Instead, we want to estimate the optimal policy $\pi^*$ based on sample transitions. 

We construct an unbiased sample of $P=\{P_a\}$ by generating one sample transition from every $(s,a)$, i.e., $\xi\!:=\!\big\{(s,a,s',\hat r_{sa})\!:\! \forall s\!\in\!\cS, a\!\in\!\cA, s'\!\sim\! P_a(s,\cdot)\big\}.$ In other words, each $\xi$ consists of $|\cS||\cA|$ sample transitions. 
Thus we obtain an sample transition matrix $P_{\xi }=\{P_a\}$ where $P_{\xi,a}(s,s') = 1$ if $s'$ is sampled and $P_{\xi,a}(s,s') = 0$ otherwise. Thus, we can define a stochastic sample of the objective function of \eqref{prob:aMDP} as 
\begin{equation}
\label{defn:aMDP-sto-fun-1}
\Phi_{\xi}(x,y) = \sum_{s\in\cS}\sum_{a\in\cA}y_{sa}\hat r_{sa} - \sum_{a\in\cA}y_a^T(P_{\xi,a}-I)x.
\end{equation}
It is easy to see that $\Phi(x,y) = \E_{\xi}\left[\Phi_{\xi}(x,y)\right]$.

\subsection{Efficiency of the Empirical Optimal Policy}
To handle the bilinear objective function, we consider the regularized empirical saddle point (R-ESP) problem, given by 
\begin{equation}
\label{prob:ESP-aMDP-reg}
\min_{x\in\cX}\max_{y\in\cY} \frac{\alpha_x}{2}\|x\|_2^2 + \Phi_{\Gamma}(x,y) - \alpha_y\sum_{s,a}y_{sa}\log(y_{sa}),
\end{equation}
where $\Phi_{\Gamma}(x,y):= \frac{1}{|\Gamma|}\sum_{\xi\in\Gamma}\Phi_\xi(x,y)$ is the empirical objective, $\alpha_x,\alpha_y$ are to be chosen later. Let  $(\bar x,\bar y)$ be the solution to the R-ESP problem \eqref{prob:ESP-aMDP-reg}. Then we obtain the empirical optimal policy 
$\bar \pi$ given by 
$$\bar{\pi}(a|s): = \bar y_{sa}/\big(\sum_{a'\in\cA}\bar y_{sa'}\big),\qquad \forall s\in\cS,a\in\cA.$$

The entropy regularizer $\sum_{s,a}y_{sa}\log(y_{sa})$ plays an important role in the analysis. It is $1$-strongly convex in $L_1$-norm due to the Pinsker’s inequality. To analyze the efficiency of $\bar\pi$, we will apply the generalization theory for SSP problem by choosing the norms as $\|\cdot\|:=\|\cdot\|_2$ and $\normm{\cdot}:=\|\cdot\|_1$, respectively.   

\begin{theorem}[Sample Efficiency of Empirical Optimal Policy]
	\label{theorem:complexity-aMDP}
	Let  $\alpha_x = \frac{\tau^{3/2}}{\sqrt{n}|\cS|t_{mix}}$,  $\alpha_y = \frac{t_{mix}}{\sqrt{n\log(|\cS||\cA|)}}$. Then the empirical optimal policy $\bar \pi$ satisfies
	$$\E\left[\hat v^*-v^{\bar \pi}\right] \leq \cO\left(\frac{t_{mix}\tau}{\sqrt{n}} \cdot\left(\tau^{1.5}+\sqrt{\log(|\cS||\cA|)}\right)\right).$$
	Consequently, to guarantee that $\E\left[\hat v^*-v^{\bar \pi}\right]\leq\epsilon$, we need $n= \Omega\left(\frac{t_{mix}^2}{\epsilon^2}\left(\tau^5 + \tau^3\log(|\cS||\cA|)\right)\right)$. Since each $\xi$ consists of $|\cS||\cA|$ samples of state transitions, the total sample complexity will be $|\cS||\cA|\cdot n = \tilde{\cO}\left(\frac{t_{mix}^2\tau^5|\cS||\cA|}{\epsilon^2}\right)$. The $|\cS||\cA|/\epsilon^{2}$ dependence in this bound is optimal.
\end{theorem}

Theorem \ref{theorem:complexity-aMDP} has several implications:
\begin{itemize}
	\item The regularized empirical optimal policy $\bar \pi$ achieves a near-optimal sample complexity, which matches known upper/lower bounds in their dependences on $|\cS|,|\cA|,\epsilon$ \cite{chen2017lower}. This result is somewhat surprising: It means that one can simply compute an empirical MDP and solve it for estimating the optimal policy. This approach is conceptually simple, yet has satisfying error bound. 
	\item Also note that the transition matrix $P$ contains $|\cS|^2|\cA|$ unknown variables, but the policy error of $\bar \pi$ scales with $|\cS||\cA|$ which is significantly smaller. Namely, one does not need to estimate the full matrix $P$ but can still get a good policy estimator by solving the R-ESP.
	\item
	The proof of Theorem \ref{theorem:complexity-aMDP} is nontrivial because we want to evaluate $v^{\bar \pi}$, which is the average reward of the state-transition process if $\bar\pi$ is implemented. The first step of the proof is to apply the result of Lemma \ref{theorem:weak-gen-reg} to the R-ESP \eqref{prob:ESP-aMDP-reg}, so that we obtain a WGM upper bound as
	\begin{equation}
	\label{eqn:aMDP-WGM}
	\Delta^w(\bar x,\bar y) \leq \cO\left(t_{mix}\cdot\big(\tau^{1.5}+\sqrt{\log(|\cS||\cA|)}\big)\cdot n^{-\frac{1}{2}}\right)
	\end{equation} 
	Then, we exploit the stationarity condition of the MDP \eqref{prob:aMDP} and prove that 
	\begin{equation}
	\label{eqn:aMDP-complementarity}
	\hat v^* - \E\Big[\sum_{a\in\cA}\bar{y}_a^\top(( P_a-I)x^*+r_a)\Big]  = \E\big[\Phi(\bar{x},y^*) - \Phi(x^*,\bar y)\big] \leq \Delta^w(\bar x,\bar y).
	\end{equation}
	Finally, we use the uniform ergodicity property of the MDP to show that 
	$\E\left[\hat v^*-v^{\bar \pi}\right] \leq \tau\E\big[\hat v^* - \sum_{a\in\cA}\bar{y}_a^\top((P_a-I)x^*+r_a)\big]$ (\cite{chen2018scalable}), which further leads to our theorem. Proofs of \eqref{eqn:aMDP-WGM}, \eqref{eqn:aMDP-complementarity} are given in Appdices \ref{appdx:Weak-gen-aMDP},\ref{appdx:complemetarity-aMDP}.
\end{itemize}


\section{Application to Stochastic Games}
\label{sec:Game}
Consider the two-player stochastic matrix game problem \eqref{prob:game}:
\begin{equation*}
\min_{x\in\Delta_{N_1}}\max_{y\in\Delta_{N_2}} x^\top\E_{\xi}[A_\xi]y\qquad\mbox{with}\qquad\Delta_{N_i}:=\{z\in\R^{N_i}: z\geq0, \mathbf{1}^\top z=1\},~i=1,2,
\end{equation*}
where $x,y$ denote the mixed strategies of players 1 and 2, respectively. Based on $n$ i.i.d. samples of the payoff matrix (with $nN_1N_2$ individual sample payoffs), we estimate the Nash equilibrium $(x^*,y^*)$ by constructing the following R-ESP problem:
$$\min_{x\in\Delta_{N_1}}\max_{y\in\Delta_{N_2}} \frac{\sum_{i}x_i\log x_i}{\sqrt{n\log N_1}} + x^\top {\footnotesize\Big(\frac{1}{n}\sum_{i=1}^{n}A_{\xi_i}\Big)} y - \frac{\sum_{j}y_j\log y_j}{\sqrt{n\log N_2}}.$$ Let $(\bar x,\bar y)$ be the solution to the preceding R-ESP, which is referred to as the {\it empirical Nash equilibrium}. Then the following theorem holds. 
\begin{theorem}
	\label{theorem:game}
	Assume $\max_{i,j}|A_\xi(i,j)|\leq 1$ almost surely. 
	Therefore, 
	$$ \E[ x^\top A \bar y ] - \cO\left(\sqrt{\log(N_1N_2)/n}\right) 
	\leq \E[\bar x^\top A\bar y] \leq \E[\bar x^\top Ay ] + \cO\left(\sqrt{\log(N_1N_2)/n}\right),$$ 
	for any $x\in\Delta_{N_1}$ and $y\in\Delta_{N_2}$.
\end{theorem}
The theorem means that the empirical strategy $(\bar x,\bar y)$ is an $\epsilon$-Nash equilibrium with high probability, as long as the total number of sample payoffs is greater than $\tilde{\cO}(N_1N_2\epsilon^{-2})$. This sample complexity is statistically optimal.

\bibliography{gen_spp}
\bibliographystyle{abbrv}


\appendix
$~$\,\,\,\,\,\,\,\,\,\,\,\,\,\,\,\,\,\,\,\,\,\,\,\,\,\,\,\,\,\,\,\,\,\,\,\,\,\,\,\,\,\,\,\,\,\,\,\,\,\,\,\,\,\,\,\,\,\,\,\,\,\,\,\,\,\,\,\,\,\,\,\,\,\,\,\,\,\,\,\,\,\,\,\,\,\,\,\,\,\,\,\,\,\,\,\,\,\,\,\,\,\,\,\,{\Large\textbf{Appendix}}
\section{Proof of Section \ref{sec:main}}
\subsection{Proof of Theorem \ref{theorem:gen-SC-SC-weak}}
\label{appdx:thm-gen-sc-sc-weak}
\begin{proof}
	For the first half of the result, we have 
	\begin{eqnarray*}
		\Delta^w(\xs,\ys) 
		& = & \max_{y\in\cY} \E\left[\Phi(\xs,y)\right] - \min_{x\in\cX} \E\left[\Phi(x,\ys)\right] \\
		& \geq & \E\left[\Phi(\xs,y^*) - \Phi(x^*,\ys)\right] \\
		& = &  \E\left[\Phi(\xs,y^*) - \Phi(x^*,y^*) + \Phi(x^*,y^*) - \Phi(x^*,\ys)\right]\\
		& \overset{\text{(a)}}{\geq} & \E\left[\frac{\mu_x}{2}\|\xs-x^*\|^2 + \frac{\mu_y}{2}\normm{\ys-y^*}^2\right]\\
		& \geq & \frac{\min\{\mu_x,\mu_y\}}{2} \cdot \E\left[d^2(\xs,\ys)\right],
	\end{eqnarray*}
	where (a) is due to the SC-SC property in Assumption \ref{assumption:convex-concave}. Rearranging the terms yields 
	$$\E\left[d^2(\xs,\ys)\right]\leq \frac{\Delta^w(\xs,\ys)}{\min\{\mu_x,\mu_y\}}.$$
	
	Next, for the proof of the weak generalization measure, we will refer to the proof of Lemma \ref{theorem:weak-gen-reg} presented in Appendix \ref{appdx:weak-gen-reg}, which prooves the generalization bound for the regularized ESP solution. The proof in Appendix \ref{appdx:weak-gen-reg} reduces to the proof of the current theorem by setting the regularizer $\Psi(x,y) = 0$.
\end{proof}

\subsection{Proof of Theorem \ref{theorem:gen-SC-SC-strong}}
\label{appdx:thm-gen-sc-sc-strong}
Before presenting the proof, let us first prove a lemma. 
\begin{lemma} 
	\label{lemma:solution-stability} 
	Suppose $\Phi$ satisfies Assumptions \ref{assumption:convex-concave} and \ref{assumption:grad-Lipschitz}. For any $y_1, y_2 \in\cY$, define $x^*(y_1)$ and $x^*(y_2)$ as $x^*(y_i) = \arg\min_{x\in\cX} \Phi(x,y_i)$ for $i = 1,2$. Similarly, for any $x_1,x_2\in\cX$, define $y^*(x_1)$ and $y^*(x_2)$ as 
	$y^*(x_i) = \argmax_{y\in\cY} \Phi(x_i,y)$, $i = 1,2$. Therefore, it holds that 
	$$\|x^*(y_1) - x^*(y_2)\| \leq \frac{L_{xy}}{\mu_x}\normm{y_1-y_2} \qquad\mbox{and}\qquad \normm{y^*(x_1)-y^*(x_2)} \leq \frac{L_{xy}}{\mu_y}\|x_1-x_2\|.$$
\end{lemma}
The proof of this lemma is presented in Appendix \ref{appdx:suport-lm-Thm2}. Now we present the proof of the Theorem \ref{theorem:gen-SC-SC-strong}.
\begin{proof}
	For the ease of notation, let us denote the Lipschitz constant 
	$$\ell_x(\xi) = \sup_{y\in\cY} \ell_x(\xi,y)\qquad\mbox{and}\qquad\ell_y(\xi) = \sup_{x\in\cX} \ell_y(\xi,x).$$
	For $\forall 1\leq i\leq n$, denote 
	$$x^*\big(\ysi\big) = \arg\min_{x\in\cX} \Phi\big(x,\ysi\big)\qquad \mbox{and}  \qquad y^*\big(\xsi\big) = \argmax_{y\in\cY} \Phi\big(\xsi,y\big).$$
	Similarly we can define $x^*\big(\ys\big)$ and $y^*\big(\xs\big).$ Therefore, 
	\begin{eqnarray}
	\label{thm:strong-1}
	& &\!\!\!\! \Phi_{\xi_i}\big(\xsi,y^*\big(\xsi\big)\!\big) - \Phi_{\xi_i}\big(x^*\big(\ysi\big),\ysi\big)\\ 
	& \overset{\text{(a)}}{\leq}  &\!\!\!\! \Phi_{\xi_i}\big(\xsi,y^*\big(\xs\big)\!\big)\! -\! \Phi_{\xi_i}\big(x^*\big(\ys\big),\ysi\big)\! +\! \ell_x(\xi_i)\big\|x^*\big(\ysi\big)-x^*\big(\ys\big)\big\|\!+ \! \ell_y(\xi_i)\normm{y^*\big(\xsi\big)\!-\!y^*\big(\xs\big)}\nonumber\\
	& \overset{\text{(b)}}{\leq} &\!\!\!\! \Phi_{\xi_i}\big(\xsi,y^*\big(\xs\big)\!\big) - \Phi_{\xi_i}\big(x^*\big(\ys\big),\ysi\big) + \frac{L_{xy}\ell_x(\xi_i)}{\mu_x}\normm{\ysi-\ys} + \frac{L_{xy}\ell_y(\xi_i)}{\mu_y}\big\|\xsi-\xs\big\|\nonumber\\
	& \overset{\text{(c)}}{\leq} &\!\!\!\! \Phi_{\xi_i}\!\big(\xs,\!y^*\big(\xs\big)\!\big)\!-\!\Phi_{\xi_i}\!\big(x^*\big(\ys\big),\!\ys\big) \!\! + \! \! \underbrace{\left(\!\frac{L_{xy}\ell_x(\xi_i)}{\mu_x} \! + \! \ell_y(\xi_i)\!\right)\!\normm{\ysi\!-\!\ys}\! \!+\! \! \left(\!\frac{L_{xy}\ell_y(\xi_i)}{\mu_y} \!+\! \ell_x(\xi_i)\!\right)\!\big\|\xsi\!-\!\xs\big\|}_{T(i)}. \nonumber
	\end{eqnarray}
	The steps (a) and (c) are due to the function Lipschitz property in Assumption \ref{assumption:Lipschitz}, and step (b) is due to Lemma \ref{lemma:solution-stability}. Consequently,
	\begin{eqnarray}
	\label{thm:strong-2}
	\Delta^s(\xs,\ys)& = &\E\left[\max_{y\in\cY} \Phi(\xs,y) - \min_{x\in\cX}\Phi(x,\ys)\right]\\
	& \overset{\text{(a)}}{=}  & \frac{1}{n}\sum_{i=1}^n\E\left[\Phi\big(\xsi,y^*\big(\xsi\big)\big) - \Phi\big(x^*\big(\ysi\big),\ysi\big)\right]\nonumber \\
	& \overset{\text{(b)}}{=} & \frac{1}{n}\sum_{i=1}^n\E\left[\Phi_{\xi_i}\big(\xsi,y^*\big(\xsi\big)\big) - \Phi_{\xi_i}\big(x^*\big(\ysi\big),\ysi\big)\right]\nonumber \\
	& \overset{\text{(c)}}{\leq}  & \E\left[\frac{1}{n}\sum_{i=1}^n\big(\Phi_{\xi_i}\big(\xs,y^*\big(\xs\big)\big)-\Phi_{\xi_i}\big(x^*\big(\ys\big),\ys\big)\big)\right]  + \frac{1}{n}\sum_{i=1}^n\E[T(i)]\nonumber\\
	& \overset{\text{(d)}}{\leq} & \frac{1}{n}\sum_{i=1}^n\E[T(i)]\nonumber.
	\end{eqnarray}
	The step (a) is because $(\xs,\ys)$ and $(\xsi,\ysi)$ are identically distributed. And the step (b) is because the independence between $\xi_i$ and $\Gamma(i)$, which indicates that 
	$$\E\left[\Phi_{\xi_i}\big(\xsi,y^*\big(\xsi\big)\big)\right] = \E\left[\E\big[\Phi_{\xi_i}\big(\xsi,y^*\big(\xsi\big)\big)\big|\xi_i\big]\right] = \E\left[\Phi\big(\xsi,y^*\big(\xsi\big)\big)\right].$$
	The independence here is a crucial point and need to be carefully handled. The step (c) is due to \eqref{thm:strong-1}. And the step (d) is because $(\xs,\ys)$ solves the ESP problem \eqref{prob:ESP}, which implies 
	$\hat\Phi_n(\xs,y) - \hat\Phi_n(x,\ys) \leq 0$ for $\forall x\in\cX,y\in\cY.$ Consequently $$\E\Big[\frac{1}{n}\sum_{i=1}^n\big(\Phi_{\xi_i}\big(\xs,y^*\big(\xs\big)\big)-\Phi_{\xi_i}\big(x^*\big(\ys\big),\ys\big)\big)\Big] = \E\left[\hat\Phi_n\left(\xs,y^*\big(\xs\big)\right) - \hat\Phi_n\left(x^*\big(\ys\big),\ys\right)\right]\leq 0.$$
	Therefore, the last step to bound $\Delta^s(\xs,\ys)$ remains as follows,
	\begin{eqnarray}
	&&\frac{1}{n}\sum_{i=1}^n\E[T(i)]\nonumber\\
	& = & \E\left[ \left(\frac{L_{xy}\ell_x(\xi_i)}{\mu_x}  +  \ell_y(\xi_i)\right)\normm{\ysi-\ys} + \left(\frac{L_{xy}\ell_y(\xi_i)}{\mu_y} + \ell_x(\xi_i)\right)\big\|\xsi-\xs\big\| \right] \nonumber\\
	& \overset{\text{(a)}}{\leq} & \E\left[\sqrt{\left(\frac{L_{xy}\ell_y(\xi_i)}{\sqrt{\mu_x}\mu_y} + \frac{\ell_x(\xi_i)}{\sqrt{\mu_x}}\right)^2+ \left(\frac{L_{xy}\ell_x(\xi_i)}{\mu_x\sqrt{\mu_y}}  +  \frac{\ell_y(\xi_i)}{\sqrt{\mu_y}}\right)^2}\cdot\sqrt{\mu_x\|\xs-\xsi\|^2 +\mu_y\normm{\ys-\ysi}^2}\right]\nonumber\\
	& \overset{\text{(b)}}{\leq} &  \sqrt{\E\left[\left(\frac{L_{xy}\ell_y(\xi_i)}{\sqrt{\mu_x}\mu_y} + \frac{\ell_x(\xi_i)}{\sqrt{\mu_x}}\right)^2+ \left(\frac{L_{xy}\ell_x(\xi_i)}{\mu_x\sqrt{\mu_y}}  +  \frac{\ell_y(\xi_i)}{\sqrt{\mu_y}}\right)^2\right]}\cdot\sqrt{\E\left[\mu_x\|\xs-\xsi\|^2 +\mu_y\normm{\ys-\ysi}^2\right]}\nonumber\\
	& \overset{\text{(c)}}{\leq} &  \sqrt{\frac{2L^2_{xy}(\ell_y^s)^2}{\mu_x\mu_y^2} + \frac{2(\ell_x^s)^2}{\mu_x}+ \frac{2L_{xy}^2(\ell_x^s)^2}{\mu_x^2\mu_y}  +  \frac{2(\ell_y^s)^2}{\mu_y}}\cdot\frac{2}{n}\cdot\sqrt{\frac{(\ell_x^s)^2}{\mu_x} + \frac{(\ell_y^s)^2}{\mu_y}}\nonumber\\
	& \leq & \frac{2\sqrt{2}}{n}\cdot\sqrt{\frac{L^2_{xy}}{\mu_x\mu_y}+1}\cdot\left(\frac{(\ell_x^s)^2}{\mu_x}+\frac{(\ell_y^s)^2}{\mu_y}\right).\nonumber
	\end{eqnarray}
	The step (a) here is due to the Cauchy-Schwartz inequality, for any two vectors $a$ and $b$, $a^\top b\leq \|a\|_2\cdot\|b\|_2$. The step (b) is the expectation version of Cauchy-Shwartz inequality, for any two random variables $a$ and $b$, $\E[ab]\leq \sqrt{\E[a^2]}\cdot\sqrt{\E[b^2]}$. And the step (c) is due to the fact that $(a+b)^2\leq 2(a^2+b^2)$ and the stability argument of Lemma \ref{lemma:stability-SC-SC}.
	
	Finally, substituting this bound into the inequality \eqref{thm:strong-2} proves the theorem. 
\end{proof}

\subsection{Proof of Theorem \ref{theorem:gen-unbounded}}
To prove the Theorem \ref{theorem:gen-unbounded}, let us first present some definition and lemmas. We define the primal function $f(x)$ and dual function $g(y)$ as well as their empirical version $f_\Gamma(x)$ and $g_\Gamma(y)$:
\begin{equation}
\label{defn:primal-dual-func}
\begin{cases}
f(x) = \max_{y} \Phi(x,y),\\
\hat f_n(x) = \max_{y}\hat\Phi_n(x,y),
\end{cases}
\quad\mbox{and}\qquad\,\,\,
\begin{cases}
g(y) = \min_{x} \Phi(x,y),\\
\hat g_n(y) = \min_{x}\hat\Phi_n(x,y).
\end{cases}
\end{equation}
For the ease of notation, we also denote
\begin{equation}
\label{defn:SunYimiao233}
x^*_n(y) = \argmin_{x} \hat\Phi_n(x,y)
\qquad\mbox{and}\qquad
y^*_n(x) = \argmax_{y}\hat\Phi_n(x,y).
\end{equation}
As a result the following property holds true. 
\begin{proposition}
	\label{proposition:primal-dual-function} 
	Under Assumption \ref{assumption:convex-concave} and \ref{assumption:grad-Lipschitz}, the primal function $f(x)$ and $\hat f_n$ are $\mu_x$-strongly convex; $\nabla f(x)$ is $L_f$-Lipschitz continuous, with $L_f:=L_x + L_{xy}^2/\mu_y$. Similarly, $g(y)$ and $\hat g_n$ are $\mu_y$-storngly concave;  $\nabla g(y)$ is $L_g$-Lipschitz continuous, with $L_g:=L_y + L^2_{xy}/\mu_x$.
\end{proposition} 
This proposition is a well known results, see e.g. \cite{sanjabi2018solving}.

\begin{lemma}
	\label{lemma:distance} 
	The squared distance from the emprical solution to the polulation solution is bounded as 
	\begin{equation} 
	\label{lm:distance-1}
	\|\xs-x^*\|_2^2\leq \frac{4}{\mu_x^2}\|\nabla_x \hat\Phi_n(x^*,y_n^*(x^*))\|_2^2\qquad\mbox{and}\qquad
	\|\ys-y^*\|_2^2\leq \frac{4}{\mu_y^2}\|\nabla_y \hat\Phi_n(x_n^*(y^*),y^*)\|_2^2. 
	\end{equation}
	and 
	\begin{equation}
	\label{lm:Bound-2nd-error-Moment-1}
	\begin{cases}
	\E\big[\|\nabla_x\hat\Phi_n(x^*,y_n^*(x^*))\|_2^2\big] \leq \frac{1}{n}\left(\frac{8L_{xy}^2}{\mu_y^2}\E\big[\|\nabla_y \Phi_{\xi}(x^*,y^*)\|_2^2\big] + 2\E\big[\|\nabla_x \Phi_{\xi}(x^*,y^*)\|_2^2\big]\right),\\
	\E\big[\|\nabla_y \hat\Phi_{n}(x^*_n(y^*),y^*)\|_2^2\big] \leq \frac{1}{n}\left(\frac{8L_{xy}^2}{\mu_x^2}\E\big[\|\nabla_x \Phi_{\xi}(x^*,y^*)\|^2_2\big] + 2\E\big[\|\nabla_y \Phi_{\xi}(x^*,y^*)\|_2^2\big]\right).
	\end{cases}
	\end{equation}
\end{lemma}
We provide the proof in Appendix \ref{appdx:distance}. As a result of Proposition \ref{proposition:primal-dual-function} and Lemma \ref{lemma:distance}, the proof of Theorem \ref{theorem:gen-unbounded} will follow the following argument. By the Lipschitz continuity of $\nabla f(x)$ and $\nabla g(y)$,
\begin{eqnarray*}
	& &\max_{y\in\cY}\Phi(\xs,y) - \min_{x\in\cX}\Phi(x,\ys)\\
	& = & f(\xs) - g(\ys)\\
	& \leq & f(x^*) + \frac{L_f}{2}\|\xs-x^*\|_2^2 - \Big(g(y^*) - \frac{L_g}{2}\|\ys-y^*\|_2^2\Big)\\
	& = & \frac{L_f}{2}\|\xs-x^*\|_2^2 + \frac{L_g}{2}\|\ys-y^*\|_2^2.
\end{eqnarray*}
Taking expectation on both sides and substituting in the values of $L_f$ and $L_g$ in Proposition \ref{proposition:primal-dual-function} proves the Theorem. 

\subsection{Proof of Lemma \ref{lemma:stability-reg}}
\label{appdx:stability-reg}
Before starting the proof, we emphasize that this is the proof of both Lemma \ref{lemma:stability-SC-SC} and Lemma \ref{lemma:stability-reg}. To get the proof of Lemma \ref{lemma:stability-SC-SC}, we can set $\Psi = 0$ and $\nu_x = \nu_y = 0$ in the following proof. 
\begin{proof}
	First, to shorten the notation in this proof, let us denote
	$$A =  \frac{\ell_x(\xi_i,\ysi) + \ell_x(\xi_i',\ys)}{n}\big\|\xs-\xsi\big\| +  \frac{\ell_y(\xi_i,\xsi) + \ell_y(\xi_i',\xs)}{n}\normm{\ys-\ysi}.$$ 
	In parallel to $\hat \Phi_n$, we define $$\hat \Phi_{n,i}(x,y) = \frac{1}{n}\sum_{\xi\in\Gamma(i)}\Phi_{\!\xi}(x,y).$$
	Then we have 
	\begin{eqnarray}
	\label{lm:stability-SC-SC-1}
	& & \hat\Phi_n\big(\xsi,\ys\big) + \Psi\big(\xsi,\ys\big) - \hat\Phi_n\big(\xs,\ysi\big) -\Psi\big(\xs,\ysi\big) \\
	& = & \frac{1}{n} \sum_{j=1}^n \Big(\Phi_{\xi_j}\big(\xsi,\ys\big)-\Phi_{\xi_j}\big(\xs,\ysi\big)\Big)+ \Psi\big(\xsi,\ys\big)-\Psi\big(\xs,\ysi\big) \nonumber\\
	& = & \frac{1}{n} \Big(\Phi_{\xi_i'}\big(\xsi,\ys\big)-\Phi_{\xi_i'}\big(\xs,\ysi\big)+\sum_{j=1,j\neq i}^n \big(\Phi_{\xi_j}\big(\xsi,\ys\big)-\Phi_{\xi_j}\big(\xs,\ysi\big)\big)\Big)\nonumber\\
	& &  + \frac{1}{n}\Big(\Phi_{\xi_i}\big(\xsi,\ys\big)-\Phi_{\xi_i}\big(\xsi,\ysi\big) + \Phi_{\xi_i}\big(\xsi,\ysi\big)-\Phi_{\xi_i}\big(\xs,\ysi\big)\Big)\nonumber\\
	& & - \frac{1}{n}\Big(\Phi_{\xi_i'}\big(\xsi,\ys\big)-\Phi_{\xi_i'}\big(\xs,\ys\big) +\Phi_{\xi_i'}\big(\xs,\ys\big) - \Phi_{\xi_i'}\big(\xs,\ysi\big)\Big)\nonumber\\
	&& + \Psi\big(\xsi,\ys\big)-\Psi\big(\xs,\ysi\big)\nonumber\\
	&\overset{\text{(a)}}{\leq} & \Big(\hat\Phi_{n,i}\big(\xsi,\ys\big)+ \Psi\big(\xsi,\ys\big) - \hat\Phi_{n,i}\big(\xs,\ysi\big)-\Psi\big(\xs,\ysi\big)\Big) + A\nonumber\\
	& = & \Big(\hat\Phi_{n,i}\big(\xsi,\ys\big)+ \Psi(\xsi,\ys) - \hat\Phi_{n,i}\big(\xsi,\ysi\big)-\Psi\big(\xsi,\ysi\big)\Big) \nonumber\\
	& & + \Big(\hat\Phi_{n,i}\big(\xsi,\ysi\big)+\Psi\big(\xsi,\ysi\big)- \hat\Phi_{n,i}\big(\xs,\ysi\big)+\Psi\big(\xs,\ysi\big)\Big) + A \nonumber\\
	& \overset{\text{(b)}}{\leq} & -\frac{\mu_x + \nu_x}{2}\|\xs-\xsi\|^2 - \frac{\mu_y+\nu_y}{2}\normm{\ys-\ysi}^2+ A.\nonumber
	\end{eqnarray}
	The step (a) is due to Lipschitz continuity of $\Phi_{\xi_i}, \Phi_{\xi'_i}$. The step (b) is due the $(\mu_y+\nu_y)$-strong concavity of $\hat\Phi_{n,i}\big(\xsi,\cdot\big)+ \Psi\big(\xsi,\cdot\big)$ and the fact that 
	$$\ysi = \argmax_{y\in\cY} \hat\Phi_{n,i}\big(\xsi,y\big)+ \Psi\big(\xsi,y\big).$$
	Hence 
	$$\hat \Phi_{n,i}\big(\xsi,\ys\big)+ \Psi\big(\xsi,\ys\big) - \hat\Phi_{n,i}\big(\xsi,\ysi\big)-\Psi\big(\xsi,\ysi\big)\leq -\frac{\mu_y+\nu_y}{2}\normm{\ys-\ysi}^2.$$
	The other part of argument on $-\frac{\mu_x + \nu_x}{2}\|\xs-\xsi\|^2$ is similar. On the other hand, similar to the argument of step (b) above, because $\xs, \ys$ and solves the strongly convex and strongly concave R-ESP problem \eqref{prob:ESP-reg}, we also have 
	\begin{eqnarray}
	\label{lm:stability-SC-SC-2}
	&&\hat\Phi_n\big(\xsi,\ys\big)+\Psi\big(\xsi,\ys\big) - \hat\Phi_n\big(\xs,\ysi\big) -\Psi\big(\xs,\ysi\big) \nonumber\\
	&\geq&\frac{\mu_x + \nu_x}{2}\|\xs-\xsi\|^2 + \frac{\mu_y+\nu_y}{2}\normm{\ys-\ysi}^2.
	\end{eqnarray}
	Combining the \eqref{lm:stability-SC-SC-1} and \eqref{lm:stability-SC-SC-2} yields 
	\begin{eqnarray*} 
		& &\!\! \!\!(\mu_x+\nu_x)\|\xs-\xsi\|^2 + (\mu_y+\nu_y)\normm{\ys-\ysi}^2\\
		& \!\!\leq & \!\!\!\!\frac{\ell_x\big(\xi_i,\ysi\big) + \ell_x\big(\xi_i',\ys\big)}{n} \|\xs-\xsi\| +  \frac{\ell_y\big(\xi_i,\xsi\big) + \ell_y\big(\xi_i',\xs\big)}{n}\normm{\ys-\ysi}\\
		& \!\!\leq & \!\!\!\!\frac{1}{n}\sqrt{\!\frac{\big(\ell_x\!\big(\xi_i,\!\ysi\big) \!+ \ell_x\!\big(\xi_i',\!\ys\big)\!\big)^2\!}{\mu_x+\nu_x} \!+\! \frac{\!\big(\ell_y\!\big(\xi_i,\!\xsi\big) \!+\! \ell_y\!\big(\xi_i',\!\xs\big)\!\big)^2\!}{\mu_y+\nu_y}}\!\cdot\! \sqrt{(\mu_x\!+\!\nu_x)\|\xs\!-\!\xsi\|^2 \!+\! (\mu_y\!+\!\nu_y)\normm{\ys\!-\!\ysi}^2}.
	\end{eqnarray*}
	Where the last row uses the Caucy-Schwartz inequality. Dividing both sides by \\$\sqrt{(\mu_x+\nu_x)\|\xs-\xsi\|^2 + (\mu_y+\nu_y)\normm{\ys-\ysi}^2}$ proves this lemma. 
\end{proof}

\subsection{Proof of Lemma \ref{theorem:weak-gen-reg}}
\label{appdx:weak-gen-reg}
\begin{proof}
	By the function Lipschitz continuity of Assumption \ref{assumption:Lipschitz}, for any $1\leq i\leq n$, and for any $x\in\cX$ and $y\in\cY$,  
	\begin{eqnarray}
	\label{thm:gen-reg-0}
	\Phi_{\xi_i}\!\big(\xsi,y\big)\!-\!\Phi_{\xi_i}\!\big(x,\ysi\big)
	\!\leq\! 
	\Phi_{\xi_i}\!\big(\xs,y\big)\!-\!\Phi_{\xi_i}\big(x,\ys\big) \!+\! \underbrace{\ell_x\big(\xi_i,y\big)\|\xs\!-\!\xsi\| \!+\! \ell_y\big(\xi_i,x\big)\normm{\ys\!-\!\ysi}}_{T(i)}.
	\end{eqnarray}
	As a result, we have 
	\begin{eqnarray}
	\label{thm:gen-reg-1}
	& & \frac{1}{n}\sum_{i=1}^n\E\big[\Phi\big(\xsi,y\big) - \Phi\big(x,\ysi\big)\big] + \Psi\big(\xs,y\big) - \Psi\big(x,\ys\big)\\
	& \overset{\text{(a)}}{=} & \frac{1}{n}\sum_{i=1}^n\E\big[\Phi_{\xi_i}\big(\xsi,y\big) - \Phi_{\xi_i}\big(x,\ysi\big)\big] + \Psi\big(\xs,y\big) - \Psi\big(x,\ys\big)\nonumber\\
	& \overset{\text{(b)}}{\leq} & \E\left[\frac{1}{n}\sum_{i=1}^n\big(\Phi_{\xi_i}\big(\xs,y\big)-\Phi_{\xi_i}\big(x,\ys\big)\big)\right] + \Psi\big(\xs,y\big) - \Psi\big(x,\ys\big) +  \frac{1}{n}\sum_{i=1}^n\E[T(i)],\nonumber
	\end{eqnarray}
	The step (a) is due to the fact that $(\xsi,\ysi)$ is independent from $\xi_i$, and hence one can take the expectation over the $\xi_i$'s first. And the step (b) is due to \eqref{thm:gen-reg-0}. Then, because the distribution of $\big(\xsi,\ysi\big)$ are the same as that of $\big(\xs,\ys\big)$ for any $1\leq i\leq n$. Therefore, the expectation term on the LHS of \eqref{thm:gen-reg-1} can be simplified to 
	\begin{equation}
	\label{thm:gen-reg-a}
	\frac{1}{n}\sum_{i=1}^n\E\big[\Phi\big(\xsi,y\big) - \Phi\big(x,\ysi\big)\big]  = \E\big[\Phi\big(\xs,y\big) - \Phi\big(x,\ys\big)\big].
	\end{equation}
	Second, the first term on the RHS of \eqref{thm:gen-reg-1} is actually
	\begin{eqnarray}
	\label{thm:gen-reg-b}
	& & \E\left[\frac{1}{n}\sum_{i=1}^n(\Phi_{\xi_i}(\xs,y)-\Phi_{\xi_i}(x,\ys))\right] + \Psi(\xs,y) - \Psi(x,\ys)\\
	& = & \E\left[\hat\Phi_n(\xs,y)-\hat\Phi_n(x,\ys)+ \Psi(\xs,y) - \Psi(x,\ys)\right]\nonumber \\
	& \leq & 0\nonumber, 
	\end{eqnarray}
	for $\forall (x,y)\in\cX$ and $y\in\cY$, which is because $(\xs,\ys)$ solves the R-ESP problem \eqref{prob:ESP-reg}. Third, because the distributions of $T(i)$'s are the same, for the second term on the RHS of \eqref{thm:gen-reg-1}, we have 
	\begin{eqnarray}
	\label{thm:gen-reg-c}
	& & \frac{1}{n}\sum_{i=1}^n\E[T(i)] \\
	& = & \E\big[\ell_x(\xi_i,y)\|\xs-\xsi\| + \ell_y(\xi_i,x)\normm{\ys-\ysi}\big]\nonumber\\
	& \overset{\text{(a)}}{\leq} & \E\left[\sqrt{\frac{\ell^2_x(\xi_i,y)}{\mu_x+\nu_x} + \frac{\ell^2_y(\xi_i,x)}{\mu_y+\nu_y}}\cdot\sqrt{(\mu_x+\nu_x)\|\xs-\xsi\|^2 +(\mu_y+\nu_y)\normm{\ys-\ysi}^2}\right]\nonumber\\
	& \overset{\text{(b)}}{\leq} & \sqrt{\E\left[\frac{\ell^2_x(\xi_i,y)}{\mu_x+\nu_x} + \frac{\ell^2_y(\xi_i,x)}{\mu_y+\nu_y}\right]}\cdot\sqrt{\E\left[(\mu_x+\nu_x)\|\xs-\xsi\|^2 +(\mu_y+\nu_y)\normm{\ys-\ysi}^2\right]} \nonumber\\
	& \overset{\text{(c)}}{\leq} & \sqrt{\E\left[\frac{\ell^2_x(\xi_i,y)}{\mu_x+\nu_x} + \frac{\ell^2_y(\xi_i,x)}{\mu_y+\nu_y}\right]}\cdot \frac{1}{n}\sqrt{\E\left[\frac{(\ell_x(\xi_i,\ysi) + \ell_x(\xi_i',\ys))^2}{\mu_x+\nu_x} + \frac{(\ell_y(\xi_i,\xsi) + \ell_y(\xi_i',\xs))^2}{\mu_y+\nu_y}\right]}\nonumber\\
	& \overset{\text{(d)}}{\leq} & \sqrt{\frac{2(\ell^w_x)^2}{\mu_x+\nu_x} + \frac{2(\ell^w_y)^2}{\mu_y+\nu_y}}\cdot \frac{1}{n}\sqrt{\frac{2(\ell^w_x)^2 + 2(\ell^w_x)^2}{\mu_x+\nu_x} + \frac{2(\ell^w_y)^2 + 2(\ell^w_y)^2}{\mu_y+\nu_y}}\nonumber\\
	& = & \frac{2\sqrt{2}}{n}\cdot\left(\frac{(\ell^w_x)^2}{\mu_x+\nu_x} + \frac{(\ell_y^w)^2}{\mu_y+\nu_y}\right).\nonumber
	\end{eqnarray}
	The step (a) uses the vector Cauchy-Schwartz inequality $a^\top b\leq\|a\|_2\cdot\|b\|_2$ for some vectors $a$ and $b$. The step (b) uses the expectation version of Cauchy-Schwartz inequality $\E[ab]\leq \sqrt{\E[a^2]}\cdot\sqrt{\E[b^2]}$ for some random variables $a$ and $b$. The step (c) is due to Lemma \ref{lemma:stability-reg}. And the step (d) is due to Assumption \ref{assumption:Lipschitz}, and the fact that $\big(\xsi,\ysi\big)$ is independent from $\xi_i$ and $\big(\xs,\ys\big)$ is independent from $\xi_i'$. Finally, substituting \eqref{thm:gen-reg-a}, \eqref{thm:gen-reg-b}, and \eqref{thm:gen-reg-c} into \eqref{thm:gen-reg-1} provides the following result:
	\begin{eqnarray*}
		\E\big[\Phi\big(\xs,y\big)\big] - \E\big[\Phi\big(x,\ys\big)\big] +\Psi\big(\xs,y\big) - \Psi\big(x,\ys\big) \leq \frac{2\sqrt{2}}{n}\cdot\left(\frac{(\ell^w_x)^2}{\mu_x+\nu_x} + \frac{(\ell_y^w)^2}{\mu_y+\nu_y}\right), \quad\mbox{for}\quad\forall x,y.
	\end{eqnarray*}
	Due to the bound the regularizer, we know $|\Psi(\xs,y) - \Psi(x,\ys)|\leq 2R$. Note that the above inequality is true for any $x$ and $y$. Therefore, we prove the overal result that
	\begin{eqnarray*}
		\max_{y\in\cY}\E\big[\Phi\big(\xs,y\big)\big] - \min_{x\in\cX}\E\big[\Phi\big(x,\ys\big)\big] \leq \frac{2\sqrt{2}}{n}\cdot\left(\frac{(\ell^w_x)^2}{\mu_x+\nu_x} + \frac{(\ell_y^w)^2}{\mu_y+\nu_y}\right)+ 2R.
	\end{eqnarray*}
	This completes the proof. 
\end{proof}

\section{Proof of Section \ref{sec:aMDP}}
\subsection{Assumptions on fast mixing time and uniform ergodicity}
\label{appdx:aMDP-assumptions}
\begin{assumption}[Uniformly bounded ergodicity]
	\label{assumption:Ergodic} The
	Markov decision process is ergodic under any stationary
	policy $\pi$, and there exists $\tau>1$ such that
	$$\frac{1}{\sqrt{\tau}|\cS|}\cdot\mathbf{1}\leq \sum_{a\in\cA} y^\pi_a\leq \frac{\sqrt{\tau}}{|\cS|}\cdot\mathbf{1},$$
	where $y^\pi$ is the stationary state-action distribution under the policy $\pi$.
\end{assumption}
\begin{assumption}[Fast mixing time] 
	\label{assumption:mixTime}
	There exists a constant $t_{mix}$ such that for any stationary policy $\pi$, 
	$$t_{mix}\geq\min_t\big\{t:\|P_\pi^t(s,\cdot) - \sum_{a\in\cA}y^\pi_a\|_{TV}\leq 1/4, \forall s\in\cS\big\},$$
	where $\|\cdot\|_{TV}$ is the total variation norm, $P_\pi(s,s') = \sum_{a\in\cA}\pi(a|s)P_a(s,s')$ is the transition probability matrix under policy $\pi$ and $P_\pi^t(s,s')$ is the $t$-step transition probability from $s$ to $s'$.
\end{assumption}

\subsection{Proof of inequality \eqref{eqn:aMDP-WGM}}
\label{appdx:Weak-gen-aMDP}
To compute the upperbound of $\Delta^w(\bar x,\bar y)$, we will first need the following proposition on the Lipschitz constants $\ell_x^w$ and $\ell_y^w$, whose proof is delegated to Appendix \ref{appdx:aMDP-Lip}.
\begin{proposition}
	\label{proposition:paras-aMDP}
	For any $\xi$ there exist constants $\ell_x(\xi,y)$ and $\ell_y(\xi,x)$ s.t. $\Phi_\xi(\cdot,y)$ is $\ell_x(\xi,y)$-Lipschitz under $L_2$-norm, and $\Phi_\xi(x,\cdot)$ is $\ell_y(\xi,x)$-Lipschitz under $L_1$-norm. Moreover,
	$$	(\ell_x^w)^2: = \sup_{y\in\cY}\E_\xi[\ell^2_x(\xi,y)] = \cO\left(\tau^{3}/|\cS|\right)\qquad\mbox{and}\qquad(\ell_y^w)^2: = \sup_{x\in\cX}\E_\xi[\ell^2_y(\xi,x)] = \cO\big(t^2_{mix}\big). $$ 
\end{proposition}

For the rest of the proof, it suffices to specify the following details for Lemma \ref{theorem:weak-gen-reg}. For the $\Phi_{\xi}$'s, $\mu_x = \mu_y = 0$. The norm $\|\cdot\|$ is the $L_2$-norm $\|\cdot\|_2$ and the norm $\normm{\cdot}$ os tje $L_1$-norm $\|\cdot\|_1$. We set the regularizer to be 
$$\Psi(x,y) = \frac{\alpha_x}{2}\|x\|_2^2 - \alpha_y\sum_{sa}y_{sa}\log y_{sa}.$$
$\Psi(\cdot,y)$ is $\nu_x$-strongly convex in $x$ under the norm $\|\cdot\|_2$ with $\nu_x = \alpha_x$. $\Psi(x,\cdot)$ is $\nu_y$-strongly concave in $y$ under the norm $\|\cdot\|_1$ with $\nu_y = \alpha_y$. Furthermore, for any $x\in\cX$ and $y\in\cY$, we know 
\begin{eqnarray*}
	R & = &\max_{x\in\cX}\max_{y\in\cY}|\Psi(x,y)|\\
	&\leq& \max_{x\in\cX}\frac{\alpha_x}{2}\|x\|_2^2 + \max_{y\in\cY}\big| \alpha_y\sum_{sa}y_{sa}\log y_{sa}\big| \\
	& = & 2\alpha_x|\cS|t_{mix}^2 + \alpha_y\log (|\cS||\cA|).
\end{eqnarray*}

From Proposition \ref{proposition:paras-aMDP}, we know $(\ell^w_x)^2 = \cO(\tau^3/|\cS|)$ and $(\ell^w_y)^2 = \cO(t^2_{mix})$. We get 
$$\max_{y\in\cY}\E[\Phi(\xs,y)] - \min_{x\in\cX}\E[\Phi(x,\ys)] \leq \cO\left(\frac{\tau^3}{n|\cS|\alpha_x} + \frac{t_{mix}^2}{n\alpha_y}+ \alpha_x|\cS|t_{mix}^2  + \alpha_y\log(|\cS||\cA|)\right).$$
To minimize the RHS, we  set 
$$\alpha_x = \frac{\tau^{3/2}}{\sqrt{n}|\cS|t_{mix}}\qquad\mbox{and}\qquad\alpha_y = \frac{t_{mix}}{\sqrt{n\log(|\cS||\cA|)}}.$$
This imediately yields  
$$\Delta^w(\xs,\ys) = \max_{y\in\cY}\E[\Phi(\xs,y)] - \min_{x\in\cX}\E[\Phi(x,\ys)] \leq \cO\left(\frac{t_{mix}}{\sqrt{n}} \left(\tau^{1.5} + \sqrt{\log(|\cS||\cA|)}\right)\right).$$

\subsection{Proof of Inequality \eqref{eqn:aMDP-complementarity}}
\label{appdx:complemetarity-aMDP}
To prove this result, let us first introduce the primal and dual linear programming formulations of the aMDP problem, which are 
\begin{eqnarray*}
	\label{prob:aMDP-primal}
	\text{(Primal-LP)}\qquad\qquad\min_{\hat v\in\R, x\in\R^{|\cS|}} \hat v\quad
	\mathrm{s.t.}  \quad\hat v\cdot\mathbf{1}  + (I-P_a)x - r_a\geq0 \quad\mbox{ for }\quad \forall a\in\cA.\qquad\qquad
\end{eqnarray*}
and 
\begin{eqnarray*}
	\text{(Dual-LP)}\qquad\qquad\min_{y\in\R^{|\cS|\times|\cA|}}\quad  \langle y, r \rangle,\quad\mathrm{s.t.}\quad y\geq0, \quad \|y\|_1 = 1,\,\, \sum_{a\in\cA} (I-P_a^\top)y_a = 0.\qquad\quad
\end{eqnarray*}
Then our saddle point problem \eqref{prob:aMDP} is the min-max formulation of this primal-dual LP pair. Let $(\hat v^*, x^*)$ be the optimal solution  to the (Primal-LP) \eqref{prob:aMDP-primal} and let $y^*$ be the optimal soltuion to the (Dual-LP), then $(x^*,y^*)$ forms the saddle point of our problem \eqref{prob:aMDP}. The following set of conditions are satisfied
\begin{equation}
\begin{cases}
\hat v^*  + (I-P_a)x^* - r_a\geq0 \quad\mbox{ for }\quad \forall a\in\cA, & \mbox{(Primal feasibility)}\\
y^*\geq0, \quad \|y^*\|_1 = 1, \sum_{a\in\cA} (I-P_a^\top)y^*_a = 0, & \mbox{(Dual feasibility)}\\
\langle y^*,r\rangle = \hat v^*,& \mbox{(Complementarity slackness)}
\end{cases}
\end{equation}
With these preliminary results, let us now provide the proof of this lemma. 
\begin{proof}
	Note that $\Phi(x,y): =\ip{y,r}+\sum_{a\in\cA} y_a^\top(P_a-I)x$, by direct computation, we have 
	\begin{eqnarray}
	& & \max_{y\in\cY}\E\left[\Phi(\bar x,y)\right] - \min_{x\in\cX}\E\left[\Phi(x,\bar{y})\right] \nonumber\\
	& \geq & \E\left[\Phi(\bar x, y^*) - \Phi(x^*,\bar{y})\right] \nonumber\\
	& = & \E\left[\sum_{a\in\cA}(y_a^*)^\top((P_a-I)\bar x+r_a) - \sum_{a\in\cA}(\bar y_a)^\top((P_a-I)x^*+r_a) \right]\nonumber\\
	&  \overset{\text{(a)}}{=} & \langle y^*,r\rangle- \E\left[\sum_{a\in\cA}(\bar y_a)^\top((P_a-I)x^*+r_a) \right]\nonumber\\
	&  \overset{\text{(b)}}{=} & \hat v^*- \E\left[\sum_{a\in\cA}(\bar y_a)^\top((P_a-I)x^*+r_a) \right].\nonumber
	\end{eqnarray}
	In the step (a), we apply the feasibility of $y^*$: $\sum_{a\in\cA}(y_a^*)^\top(P_a-I) = 0$. In the step (b), we applied the fact that $\hat v^* = \langle y^*,r\rangle$. Hence we complete the proof. 
\end{proof}

\section{Proof of Section \ref{sec:Game}}
\label{appdx:game}
\subsection{Proof of Theorem \ref{theorem:game}}
\begin{proof}
	The result of this theorem is a direct corollary of Lemma \ref{theorem:weak-gen-reg}. First, we will need to figure the corresponding algorithmic constants. Define
	$$\Psi(x,y) = \frac{1}{\sqrt{n\log N_1}}\sum_{i}x_i\log x_i - \frac{1}{\sqrt{n\log N_2}}\sum_j y_j\log y_j$$ as the regularizer. Then we know $\Psi$ is $\nu_x$-strongly convex under $L_1$-norm and $\nu_y$-strongly concave under $L_1$-norm, with $\nu_x = \frac{1}{\sqrt{n\log N_1}}$ and $\nu_y = \frac{1}{\sqrt{n\log N_2}}$. Morever, the magnitude of $\Psi$ is upper bounder by 
	\begin{eqnarray*}
		R & = &  \max_{x\in\Delta_{N_1}}\max_{y\in\Delta_{N_2}}|\Psi(x,y)|\\
		& = &\frac{1}{\sqrt{n\log N_1}}\log N_1 + \frac{1}{\sqrt{n\log N_2}}\log N_2\\
		& = & \frac{1}{\sqrt{n}}\left(\sqrt{\log N_1} + \sqrt{\log N_2}\right).
	\end{eqnarray*}
	Denote $\Phi_\xi(x,y) = x^\top A_\xi y$, then clearly, $\Phi_{\xi}$ is not an SC-SC function. Hence $\mu_x = \mu_y = 0$. Let the Lipschitz constants $\ell_x(\xi,y)$ and $\ell_y(\xi,x)$ be measured under the $L_1$-norm, then 
	\begin{eqnarray*}
		\ell_x(\xi,y) & = & \max_{x\in\Delta_{N_1}} \|\nabla_x \Phi_{\xi}(x,y)\|_\infty\\
		& = &\max_{x\in\Delta_{N_1}} \|A_\xi y\|_\infty\\
		&= & \max_{ij} |A_\xi(i,j)|\\
		& = &1. 
	\end{eqnarray*}
	Consequently, we have $\ell_x^w = 1$. Similarly, we have $\ell_y(\xi,x)\leq 1$ and $\ell_y^w\leq 1$. Consequently, by Lemma \ref{theorem:weak-gen-reg}, we have 
	\begin{eqnarray*}
		\Delta^w(\xs,\ys) & \leq & \frac{2\big(\ell_x^w\big)^2}{n(\mu_x+\nu_x)} + \frac{2\big(\ell_y^w\big)^2}{n(\mu_y+\nu_y)} + 2R \\
		& = & 8\frac{\sqrt{\log N_1} + \sqrt{\log N_2}}{\sqrt{n}} \\
		& \leq & 16\frac{\sqrt{\log (N_1N_2)}}{\sqrt{n}}.
	\end{eqnarray*}
	The last inequality is due to $\sqrt{\log N_1} + \sqrt{\log N_2} \leq 2\sqrt{\log N_1+\log N_2} = 2\sqrt{\log(N_1N_2)}$. Due to  the definition of $\Delta^w(\xs,\ys)$, we know that for any $x\in\Delta_{N_1}$ and $y\in\Delta_{N_2}$, 
	\begin{eqnarray*}
		\underbrace{\E\big[\bar x^\top Ay - \bar x^\top A\bar y\big]}_{\geq0} - \underbrace{\E\big[x^\top A\bar y-\bar x^\top A\bar y\big]}_{\leq0}  =  \E\big[\bar x^\top Ay  - x^\top A\bar y\big] \leq \Delta^w(\xs,\ys) \leq \cO\left(\frac{\sqrt{\log(N_1N_2)}}{\sqrt{n}}\right).
	\end{eqnarray*}
	Consequently, 
	$\E\big[\bar x^\top Ay - \bar x^\top A\bar y\big]\leq\cO\left(\frac{\sqrt{\log(N_1N_2)}}{\sqrt{n}}\right)$ for any $y$. Meaning that when the player 1 plays $\bar x$, in expectation, it does not gain much benefit for player 2 if he switches to any other fixed strategy $y$. Symmetrically, we have $\E\big[x^\top A\bar y-\bar x^\top A\bar y\big]$, meaning that when the player 2 plays the strategy $\bar y$,  in expectation, it does not cause more lost for player 2 if player 1 switches to any other fixed strategy $x$. 
\end{proof}

\section{Other supporting lemmas}
\subsection{Proof of Lemma \ref{lemma:solution-stability}}
\label{appdx:suport-lm-Thm2}
\begin{proof}
	First, let us prove the result for $\|x^*(y_1)-x^*(y_2)\|$ and the rest of the results can be proved parallelly. By the optimality condition we have 
	$$\left\langle \nabla_x \Phi(x^*(y_1),y_1), x^*(y_2)-x^*(y_1)\right\rangle \geq 0 \qquad\mbox{and}\qquad\left\langle\nabla_x \Phi(x^*(y_2),y_2), x^*(y_1)-x^*(y_2)\right\rangle\geq0.$$
	Summing this up gives 
	$$\langle\nabla_x\Phi(x^*(y_2),y_2)-\nabla_x \Phi(x^*(y_1),y_1),  x^*(y_1)-x^*(y_2)\rangle \geq 0.$$
	By the strong convexity of $\Phi(\cdot,y)$ and the $L_{xy}$-Lipschitz continuity of $\nabla_x\Phi(x,y)$ in terms of $y$, 
	\begin{eqnarray*}
		0 & \leq & \langle\nabla_x\Phi(x^*(y_2),y_2)-\nabla_x \Phi(x^*(y_1),y_1),  x^*(y_1)-x^*(y_2)\rangle\\
		& = & \langle\nabla_x\Phi(x^*(y_2),y_2)-\nabla_x \Phi(x^*(y_1),y_2),  x^*(y_1)-x^*(y_2)\rangle\\
		& & + \langle \nabla_x \Phi(x^*(y_1),y_2)-\nabla_x \Phi(x^*(y_1),y_1) ,x^*(y_1)-x^*(y_2)\rangle\\
		& \leq & -\mu_x\|x^*(y_1)-x^*(y_2)\|^2 + \|\nabla_x \Phi(x^*(y_1),y_2)-\nabla_x \Phi(x^*(y_1),y_1)\|_*\|x^*(y_1)-x^*(y_2)\|\\
		&\leq& -\mu_x\|x^*(y_1)-x^*(y_2)\|^2 + L_{xy}\normm{y_1-y_2}\cdot\|x^*(y_1)-x^*(y_2)\|.
	\end{eqnarray*}
	Consequently, 
	$$\|x^*(y_1)-x^*(y_2)\|\leq \frac{L_{xy}}{\mu_x}\normm{y_1-y_2}.$$
	The other part of the result follows the same line of proof. 
\end{proof}

\subsection{Proof of Lemma \ref{lemma:distance}}
\label{appdx:distance}
\begin{proof}
	Because $\xs = \argmin_x \hat f_n(x)$, and $\hat f_{n}(\cdot)$ is $\mu_x$-strongly convex, we have
	\begin{eqnarray*}
		0 &\geq &\hat f_{n}(\xs) - \hat f_{n}(x^*)\\
		& \geq & \langle\nabla \hat f_{n}(x^*),\xs - x^*\rangle + \frac{\mu_x}{2}\|\xs-x^*\|^2_2\\
		& = &  \langle\nabla_x\hat\Phi_n(x^*,y^*_n(x^*)),\xs - x^*\rangle + \frac{\mu_x}{2}\|\xs-x^*\|_2^2,
	\end{eqnarray*}
	where the last row is due to the Danskin's theorem. By rearranging the terms, we get
	\begin{eqnarray*}
		\frac{\mu_x}{2}\|\xs-x^*\|_2^2 \leq  -\langle\nabla_x\hat\Phi_n(x^*,y^*_n(x^*)) ,\xs - x^*\rangle
		\leq \|\nabla_x\hat\Phi_n(x^*,y^*_n(x^*))\|_2\cdot\|\xs - x^*\|_2.
	\end{eqnarray*} 
	Deviding both sides by $\|\xs - x^*\|_2$ and then square both sides proves the first inequality of \eqref{lm:distance-1}.  The second inequality for \eqref{lm:distance-1} can be proved similarly.
	
	Next, let us focus on the first inequality of \eqref{lm:Bound-2nd-error-Moment-1}.
	\begin{eqnarray}
	\label{lm:Bound-2nd-error-Moment-2}
	& & \E\big[\|\nabla_x\hat\Phi_n(x^*,y^*_n(x^*))\|_2^2\big]\\
	& = & \E\big[\| \nabla_x\hat\Phi_n(x^*,y^*_n(x^*)) - \nabla_x \hat\Phi_n(x^*,y^*) + \nabla_x \hat\Phi_n(x^*,y^*)\|_2^2\big]\nonumber\\
	& \leq & 2\E\big[\|\nabla_x\hat\Phi_n(x^*,y^*_n(x^*)) - \nabla_x \hat\Phi_n(x^*,y^*)\|_2^2\big]  + 2\E\big[\|\nabla_x \hat\Phi_n(x^*,y^*)\|_2^2\big]\nonumber\\
	& \leq & 2L_{xy}^2\E\big[\|y^*-y^*_n(x^*)\|_2^2\big] + 2\E\big[\|\nabla_x \hat\Phi_n(x^*,y^*)\|_2^2\big].\nonumber
	\end{eqnarray}
	Note that $\E\big[\nabla_x \hat\Phi_n(x^*,y^*)\big] = \nabla_x\Phi(x^*,y^*) = 0 $, we have 
	\begin{eqnarray}
	\label{lm:Bound-2nd-error-Moment-3}
	\E\big[\|\nabla_x \hat\Phi_n(x^*,y^*)\|_2^2\big] & = & \E\left[\Big\|\frac{1}{n}\sum_{\xi\in\Gamma}\nabla_x \Phi_{\xi}(x^*,y^*)-\E\left[\nabla_x\Phi_\xi(x^*,y^*)\right]\Big\|_2^2\right]\\
	& = & \frac{1}{n}\E_\xi\big[\|\nabla_x \Phi_{\xi}(x^*,y^*)\|_2^2\big].\nonumber
	\end{eqnarray}
	It is important that $L_2$-norm is used here so that the above variance equation chain holds. If another norm is used, \eqref{lm:Bound-2nd-error-Moment-3} may not be true. For example if $L_1$-norm is used, an extra multiplicative factor of  dimension will come into the bound. For the other term, note that 
	$$y^*_n(x^*) = \argmax_y \hat\Phi_n(x^*,y)\qquad\mbox{and}\qquad
	y^* = \argmax_y \Phi(x^*,y).$$
	As a result, we have 
	\begin{eqnarray*}
		0 \leq \hat\Phi_n(x^*,y^*_n(x^*)) -\hat\Phi_n(x^*,y^*)
		\leq  \langle\nabla_y \hat\Phi_n(x^*,y^*),y^*_n(x^*) - y^*\rangle - \frac{\mu_y}{2}\|y^*_n(x^*)-y^*\|_2^2.
	\end{eqnarray*}
	With slight rearranging and apply Cauchy-Schwartz inequality, we have 
	$$\|y^*_n(x^*)-y^*\|_2\leq \frac{2}{\mu_y}\|\nabla_y \hat\Phi_n(x^*,y^*)\|_2.$$
	Taking expectation on both sides and we get 
	\begin{eqnarray}
	\label{lm:Bound-2nd-error-Moment-4}
	\E\big[\|y^*_n(x^*)-y^*\|^2_2\big] & \leq & \frac{4}{\mu_y^2}\E\big[\|\nabla_y \hat\Phi_n(x^*,y^*)\|^2_2\big]\\
	& = &\frac{4}{\mu_y^2}\E\big[\|\nabla_y \hat\Phi_n(x^*,y^*) - \nabla_y\Phi(x^*,y^*)\|^2_2\big]\nonumber\\
	& \leq & \frac{4}{n\mu^2_y}\E_\xi\big[\|\nabla_y \Phi_{\xi}(x^*,y^*)\|^2_2\big].\nonumber
	\end{eqnarray}
	The argument here is parallel to that of \eqref{lm:Bound-2nd-error-Moment-3}. Combining \eqref{lm:Bound-2nd-error-Moment-2}, \eqref{lm:Bound-2nd-error-Moment-3}, and \eqref{lm:Bound-2nd-error-Moment-4}, we have 
	$$\E\big[\|\nabla_x\hat\Phi_n(x^*,y^*_n(x^*))\|^2_2\big] \leq \frac{1}{n}\left(\frac{8L_{xy}^2}{\mu_y^2}\E_\xi\big[\|\nabla_y \Phi_{\xi}(x^*,y^*)\|^2_2\big] + 2\E_\xi\big[\|\nabla_x \Phi_{\xi}(x^*,y^*)\|^2_2\big]\right).$$
	The second inequality can be proved through a completely parallel way.
\end{proof}

\subsection{Proof of Proposition \ref{proposition:paras-aMDP}}
\label{appdx:aMDP-Lip}
\begin{proof}
	First, note that the Lipschitz continuity in $x$ variable is measured under $L_2$-norm, and the Lipschitz continuity in $y$ variable is measured under the $L_1$-norm. Because the dual norms of $\|\cdot\|_2$ and $\|\cdot\|_1$ are $\|\cdot\|_2$ and $\|\cdot\|_\infty$ repectively, we have
	$$\ell^2_y(\xi,x) = \sup_{y\in\cY} \|\nabla_y \Phi_{\xi}(x,y)\|_\infty^2\qquad\mbox{and}\qquad\ell^2_x(\xi,y) = \sup_{x\in\cX} \|\nabla_y \Phi_{\xi}(x,y)\|_2^2.$$
	By direct calculation, we know 
	$$\nabla_{y_a} \Phi_{\xi}(x,y) = \hat r_a - x + u_a\quad\mbox{with}\quad u_{a,s} = \sum_{s'\in\cS}\delta_\xi(s,a,s')x_{s'},$$
	where $\hat r_a = [\hat{r}_{1,a}, \hat r_{2a},...,\hat r_{|\cS|a}]^\top$. Consequently, 
	\begin{eqnarray*}
		\|\nabla_{y} \Phi_{\xi}(x,y)\|_\infty  =  \max_{a\in\cA} \|\nabla_{y_a} \Phi_{\xi}(x,y)\|_\infty
		\leq  \max_{a\in\cA} \|\hat r_a\|_\infty+ \|x\|_\infty+ \|u_a\|_\infty\leq 1+4t_{mix}.
	\end{eqnarray*}
	As a result we have $(\ell_y^w)^2 = \cO(t^2_{mix}).$ For $\ell_x$, we first compute the  gradient as follows
	$$\nabla_{x} \Phi_{\xi}(x,y) = -\sum_{a\in\cA}y_a + w\quad\mbox{with}\quad w_{s'} = \sum_{s,a}y_{sa}\delta_\xi(s,a,s').$$
	Consequently, for any fixed $y\in\cY$,
	\begin{eqnarray}
	\label{prop:aMDP-1}
	\E_\xi[\ell^2_x(\xi,y)] & = & \|\nabla_{x} \Phi_{\xi}(x,y)\|_2^2 \\
	& \leq & 2\big\|\sum_{a\in\cA}y_a\big\|_2^2 + 2\E_{\xi}[\|w\|_2^2]\nonumber\\
	& = & 2\big\|\sum_{a\in\cA}y_a\big\|_2^2 + 2\|\E[w]\|_2^2 + 2\E_\xi[\|w-\E_\xi[w]\|_2^2].\nonumber
	\end{eqnarray}
	Note that $\sum_{a\in\cA}y_a\leq \frac{\sqrt{\tau}}{|\cS|}\cdot\mathbf{1}$, we know $\|\sum_{a\in\cA}y_a\|_2^2\leq \frac{\tau}{|\cS|}$. By directly calculating the expectation, we know 
	$$\E_{\xi}[w_{s'}] = \sum_{s,a}y_{sa}P_a(s,s').$$
	Because for a particular $s'$, $\delta_\xi(s_1,a_1,s')$ is independent from $\delta_\xi(s_2,a_2,s')$, we can compute the variance term as
	\begin{eqnarray*}
		\E_{\xi}\left[\|w-\E_\xi[w]\|_2^2\right]& = & \sum_{s'}\E_{\xi}\left[(\sum_{sa} y_{sa}(\delta_\xi(s,a,s')-P_a(s,s')))^2\right]\\
		& = & \sum_{s'}\sum_{sa} y^2_{sa}\E_{\xi}\left[(\delta_\xi(s,a,s')-P_a(s,s')))^2\right]
		\\
		& = & \sum_{s'}\sum_{sa} y^2_{sa}(1-P_a(s,s'))P_a(s,s')\\
		& \leq & \sum_{s'}\sum_{sa} y^2_{sa}P_a(s,s')\\
		& \leq & \sum_{sa} y^2_{sa}\\
		& \overset{\text{(a)}}{\leq} & \big\|\sum_a y_a\big\|_2^2\\
		& \leq & \tau/|\cS|.
	\end{eqnarray*}
	Where the step (a) is because $y\geq0$.  Now we bound the last term $\|\E_{\xi}[w]\|_2^2$. For the ease of discussion, let us define $\bar y = 0.5y + 0.5\frac{\mathbf{1}}{|\cS||\cA|}$, $\bar\lambda = \sum_{a\in\cA}\bar y_a$, and $\bar\pi(a|s) = \bar y_{sa}/\bar \lambda_s$. Similarly, we define
	$\hat y = \frac{\mathbf{1}}{|\cS||\cA|}$,	 $\hat\lambda = \frac{\mathbf{1}}{|\cS|}$, and $\hat \pi(a|s) = \frac{1}{|\cA|}$. Therefore, because both $\bar{\pi}$ and $\hat{\pi}$ are strictly positive, the corresponding Markov chains of the state transitionsis are ergodic. Hence,
	\begin{eqnarray*}
		\|\E_{\xi}[w]\|_2^2 & = & \sum_{s'}\left(\sum_{s,a}\left(y_{sa} + \frac{1}{|\cS||\cA|}-\frac{1}{|\cS||\cA|}\right)P_a(s,s')\right)^2 \\
		&\leq&2\sum_{s'}\left(\sum_{s,a}\left(y_{sa} + \frac{1}{|\cS||\cA|}\right)P_a(s,s')\right)^2 + 2\sum_{s'}\left(\sum_{s,a}\left(\frac{1}{|\cS||\cA|}\right)P_a(s,s')\right)^2\\
		&\overset{\text{(a)}}{=}&8\sum_{s'}\left(\sum_{s,a}\bar y_{sa}P_a(s,s')\right)^2 + 2\sum_{s'}\left(\sum_{s,a}\hat y_{sa}P_a(s,s')\right)^2\\
		&\overset{\text{(b)}}{=}&8\sum_{s'}\Big(\sum_{s}\bar \lambda_s\sum_a\bar\pi(a|s)P_a(s,s')\Big)^2 + 2\sum_{s'}\Big(\sum_{s}\hat\lambda_s\sum_a\hat\pi(a|s)P_a(s,s')\Big)^2\\
		& \overset{\text{(c)}}{=} & 8\sum_{s'}\Big(\sum_s \bar\lambda_sP_{\bar\pi}(s,s')\Big)^2+2\sum_{s'}\Big(\sum_s \hat\lambda_sP_{\hat\pi}(s,s')\Big)^2\\
		& = & 8\|P_{\bar\pi}^\top\bar\lambda\|_2^2 + 2\|P_{\hat\pi}^\top\hat\lambda\|_2^2\\
		& \overset{\text{(d)}}{\leq} & \frac{10\tau^3}{|\cS|}.
	\end{eqnarray*}
	The step (a) and (b)  follows directly from the definition of $\bar y$, $\bar \lambda$, $\bar \pi$ and $\hat y$, $\hat \lambda$, $\hat{\pi}$. The step (c) is we define $P_{\bar \pi}$ to be the state transition probability matrix under the policy $\bar \pi$, and $P_{\bar\pi}(s,s'):=\sum_a\bar\pi(a|s)P_a(s,s')$; Similar argument is made for $P_{\hat{\pi}}$. Finally, the step (d) is due to the following argument. Let $\lambda^{\bar\pi}$ be the stationary state distribution under the policy $\bar\pi$, then by ergodicity property (Assumption \ref{assumption:Ergodic}) we have
	$$0\leq P_{\bar\pi}^\top \mathbf{1}\leq P_{\bar\pi}^\top(\sqrt{\tau}|\cS|\lambda^{\bar\pi}) = \sqrt{\tau}|\cS|\lambda^{\bar\pi} \leq \sqrt{\tau}|\cS|\cdot\frac{\sqrt{\tau}}{|\cS|}\cdot\mathbf{1} = \tau\mathbf{1}.$$
	As a result, $0\leq P_{\bar\pi}^\top\bar\lambda\leq P_{\bar\pi}^\top\frac{\sqrt{\tau}}{|\cS|}\mathbf{1}\leq \frac{\tau^{3/2}}{|\cS|}\mathbf{1}$ and consequently $\|P_{\bar\pi}^\top \bar\lambda\|_2^2\leq \tau^3/|\cS|$. Similarly, $\|P_{\hat\pi}^\top \hat\lambda\|_2^2\leq \tau^3/|\cS|$. Substituting the following bounds 
	$$\|\sum_a y_a\|_2^2 \leq \cO(\tau/|\cS|), \quad\quad \E_{\xi}\left[\|w-\E_\xi[w]\|_2^2\right]\leq \cO(\tau/|\cS|), \quad\mbox{and}\quad \big\|\E_\xi[w]\big\|_2^2\leq \cO(\tau^3/|\cS|)$$
	into \eqref{prop:aMDP-1} proves that 
	$\E_{\xi}[\ell_x^2(\xi,y)] \leq \cO(\tau^3/|\cS|)$. Consequently, $(\ell_x^w)^2 = \sup_{y\in\cY} \E_{\xi}[\ell_x^2(\xi,y)] \leq \cO(\tau^3/|\cS|).$ This completes the proof of this proposition. 
\end{proof}

\end{document}